\documentclass[aap,MSNbibl,seceqn,dvips]{arximspdf}
\usepackage{mathrsfs}

\doi{10.1214/10-AAP690}
\volume{21}
\issue{1}
\pubyear{2011}
\firstpage{140}
\lastpage{182}

\makeatletter

\newtheorem{Theorem}{Theorem}[section]

\newtheorem{Proposition}[Theorem]{Proposition}

\newproclaim{Definition}[Theorem]{Definition}
\newproclaim{Example}{Example}[section]

\newproclaim{Remark}[Theorem]{Remark}

\def\eqref#1{(\ref{#1})}
\makeatother

\begin{document}
\begin{frontmatter}

\title{Multivariate supOU processes}
\runtitle{Multivariate supOU processes}

\begin{aug}
\author[A]{\fnms{Ole Eiler} \snm{Barndorff-Nielsen}\ead[label=e1]{oebn@imf.au.dk}}
\and
\author[B]{\fnms{Robert} \snm{Stelzer}\ead[label=e2]{rstelzer@ma.tum.de}
\ead[label=u1,url]{http://www-m4.ma.tum.de}\corref{}}
\runauthor{O. E. Barndorff-Nielsen and R. Stelzer}
\affiliation{\AA rhus University  and  Technische Universit\"{a}t M\"{u}nchen}
\address[A]{Thiele Centre\\ Department of Mathematical Sciences\\ \AA
rhus University \\Ny Munkegade\\
DK-8000 \AA rhus C
\\ Denmark\\\printead{e1}} 
\address[B]{TUM Institute for Advanced Study \\\quad and Zentrum
Mathematik\\Technische Universit\"{a}t M\"{u}nchen\\
Boltzmannstrasse 3\\ D-85747 Garching\\ Germany\\
\printead{e2}\\\printead{u1}}
\end{aug}

\received{\smonth{5} \syear{2009}}
\revised{\smonth{2} \syear{2010}}

%
\begin{abstract}
Univariate superpositions of Ornstein--Uhlenbeck-type processes (OU),
called supOU processes, provide a class of continuous time
processes capable of exhibiting long memory behavior. This paper
introduces multivariate supOU processes and gives conditions for
their existence and finiteness of moments. Moreover, the second-order
moment structure is explicitly calculated, and examples
exhibit the possibility of long-range dependence.

Our supOU processes are defined via homogeneous and factorizable L\'evy bases.
We show that the behavior of supOU processes is particularly nice when
the mean reversion
parameter is restricted to normal matrices and especially to strictly
negative definite ones.

For finite variation L\'evy bases we are able to give conditions for
supOU processes to have locally
bounded c\`adl\`ag paths of finite variation and to show an analogue of
the stochastic differential
equation of OU-type processes, which has been suggested in \cite
{barndorffnielsen01} in the univariate case.
Finally, as an important special case, we introduce
positive semi-definite supOU processes, and we discuss the relevance of
multivariate supOU processes in applications.
\end{abstract}

%
\begin{keyword}[class=AMS]
\kwd[Primary ]{60G10}
\kwd{60H20}
\kwd[; secondary ]{60E07}
\kwd{60G51}
\kwd{60G57}.
\end{keyword}
\begin{keyword}
\kwd{L\'evy bases}
\kwd{long memory}
\kwd{normal matrices}
\kwd{Ornstein--Uhlenbeck-type processes}
\kwd{positive semi-definite stochastic processes}
\kwd{second-order moment structure}
\kwd{stochastic differential equation}.
\end{keyword}

\end{frontmatter}

\section{Introduction}
L\'evy-driven Ornstein--Uhlenbeck-type processes (OU) are extensively
used in applications as elements in continuous time models for observed
time series. One area where they are often applied is mathematical
finance (see, e.g., \cite{Contetal2004}), especially in the OU-type
stochastic volatility model of \cite{Barndorffetal2001c}. An OU-type
process is given as the solution of a stochastic differential equation
of the form
%
\begin{equation}\label{eq:ousde}
dX_t=-aX_t \,dt + dL_t
\end{equation}
with $L$ being a L\'evy process (see, e.g., \cite{Sato1999} for a
comprehensive introduction) and $a\in\mathbb{R}$. Typically, one is
interested mainly in stationary solutions of \eqref{eq:ousde}.
Provided $a>0$ and $E(\ln(|L|\vee1))<\infty$, the SDE \eqref
{eq:ousde} has a unique stationary solution given by
\[
X_t=\int_{-\infty}^te^{-a(t-s)}\,dL_s.
\]

However, in many applications the dependence structure exhibited by
empirical data is found to be not in good accordance with that of
OU-type processes which have autocorrelation functions of the form
$e^{-ah}$ for positive lags $h$. In many data sets a more complex and
often a (quasi)long memory behavior of the autocorrelation function is
encountered. OU-type processes could be replaced by fractional OU-type
processes (see \cite{Marquardt2006} or \cite{Marquardt2007}, for
instance) to have long memory effects included in the model. However,
in this case many desirable properties are lost, and, in particular,
fractional OU-type processes no longer have jumps. An alternative to
obtain long memory from OU-type processes and still to have jumps is to
add up countably many independent OU-type processes, that is,
\[
X_t=\sum_{k=1}^\infty w_i\int_{-\infty}^te^{-a_i(t-s)}\,dL_{i,s}
\]
with independent identically distributed L\'evy processes $(L_i)_{i\in
\mathbb{N}}$ and appropriate $a_i>0$, $w_i>0$ with $\sum_{i=1}^\infty w_i=1$.
Intuitively we can likewise ``add'' (i.e., integrate) up independent
OU-type processes with all parameter values $a>0$ possible. The
resulting processes are called supOU processes and have been introduced
in \cite{barndorffnielsen01} where it has also been established that
they may exhibit long-range dependence. For a comprehensive treatment
regarding the theory and use of univariate supOU processes in finance
we refer to \cite{Barndorffetal2003}.

So far supOU processes have only been considered in the univariate
case. However, in many applications it is crucial to model several time
series with a joint model, and so flexible multivariate models are
important. Therefore, in this paper we introduce and study multivariate
supOU processes. Due to the appearance of matrices and the related
peculiarities our theory is not a straightforward extension of the
univariate results. Multivariate ($d$-dimensional) OU-type processes
(see, e.g., \cite{Jureketal1993} or \cite{satoyamazato1984}) are
defined as the solutions of SDEs of the form
%
\begin{equation}\label{eq:multousde}
dX_t=AX_t\,dt+dL_t
\end{equation}
with $L$ a $d$-dimensional L\'evy process and $A$ a $d\times d$ matrix.
Provided\break $E(\ln(\|L\|\vee1))<\infty$ and all eigenvalues of $A$ have
strictly negative real part, we have again a unique stationary solution
given by
\[
X_t=\int_{-\infty}^te^{A(t-s)}\,dL_s.
\]
Intuitively our multivariate supOU processes are obtained by ``adding
up'' independent OU-type processes with all possible parameters $A$;
that is, we consider all $d\times d$ matrices $A$ with eigenvalues of
strictly negative real parts. It turns out later on that the behavior
of supOU becomes more tractable when we restrict $A$ to come only from
a nice subset, like the negative definite matrices.

The remainder of this paper is structured as follows. The next section
starts with a brief overview of important notation and conventions used
in the paper and is followed in Section \ref{sec:21} by a
comprehensive introduction into L\'evy bases and the related
integration theory, which will be needed to define supOU processes. In
Section \ref{sec:3} we first define multivariate supOU processes and
provide existence conditions in Section \ref{sec:31}. Thereafter, we
discuss the existence of moments and derive the second-order structure.
For the finite variation case we show important path properties in
Section \ref{sec:33}. Besides establishing that we have c\`adl\`ag
paths of bounded variation, we give an analogue of the stochastic
differential equation~\eqref{eq:multousde} for supOU processes and its
proof. In particular, this proves a conjecture in \cite
{barndorffnielsen01}, which has not yet been shown in any
nondegenerate set-up. We conclude that section with several examples
illustrating the behavior and properties of supOU processes and showing
that they may exhibit long memory. In Section \ref{sec:4} we use our
results to define positive semi-definite supOU processes and analyze
their properties. These processes are important for applications like
stochastic volatility modeling, since they may be used to describe the
stochastic dynamics of a latent covariance matrix. Finally, this and
other possible applications of supOU processes are discussed in Section~\ref{sec:appl}.

\section{Background and preliminaries}\label{sec2}
\subsection{Notation}

We denote the set of real $m\times n$ matrices by $M_{m,n}(\mathbb
{R})$. If
$m=n$, we simply write $M_n(\mathbb{R})$ and denote the group of invertible
$n\times n$ matrices by $GL_n(\mathbb{R})$, the linear subspace of symmetric
matrices by $\mathbb{S}_n$, the (closed) positive semi-definite cone by
$\mathbb{S}_n^+$ and the open positive definite cone by $\mathbb{S}_n^{++}$
(likewise $\mathbb{S}_n^{--}$ are the strictly negative definite
matrices,
etc.). $I_n$ stands for the $n\times n$ identity matrix. The tensor
(Kronecker) product of two matrices $A,B$ is written as $A\otimes B$.
$\operatorname{vec}$ denotes the well-known vectorisation operator
that maps the
$n\times n$ matrices to $\mathbb{R}^{n^2}$ by stacking the columns of the
matrices below one another. For more information regarding the tensor
product and $\operatorname{vec}$ operator we refer to \cite{Hornetal1991}, Chapter 4. The spectrum of a matrix is denoted by $\sigma(\cdot
)$. Finally, $A^*$ is the transpose (adjoint) of a matrix $A\in
M_{m,n}(\mathbb{R})$, and $A_{ij}$ stands for the entry of $A$ in the $i$th
row and $j$th column.

Norms of vectors or matrices are denoted by $\|\cdot\|$. If the norm
is not further specified, then it is understood that we take the
Euclidean norm or its induced operator norm, respectively. However, due
to the equivalence of all norms none of our results really depends on
the choice of norms.

For a complex number $z$ we denote by $\Re(z)$ its real part.
Moreover, the indicator function of a set $A$ is written $1_A$.

A mapping $f\dvtx V\to W$ is said to be $\mathscr{V}$--$\mathscr
{W}$-measurable if it is measurable when the $\sigma$-algebra
$\mathscr{V}$ is used on the domain $V$, and the $\sigma$-algebra
$\mathscr{W}$ is used on the range $W$. The Borel $\sigma$-algebras
are denoted by $\mathscr{B}(\cdot)$ and $\lambda$ typically stands
for the Lebesgue measure which in vector or matrix spaces is understood
to be defined as the product of the coordinatewise Lebesgue measures.

Throughout we assume that all random variables and processes are
defined on a given complete probability space $(\Omega,\mathscr
{F},P)$ equipped with an appropriate filtration when relevant.

Furthermore, we employ an intuitive notation with respect to the
(stochastic) integration with matrix-valued integrators referring to
any of the standard texts (e.g., \cite{Protter2004}) for a
comprehensive treatment of the theory of stochastic integration. Let
$(A_t)_{t\in\mathbb{R}^+}$ in $M_{m,n}(\mathbb{R})$, $(B_t)_{t\in
\mathbb{R}^+}$ in
$M_{r,s}(\mathbb{R})$ be c\`adl\`ag and adapted processes and $
(L_t)_{t\in
\mathbb{R}^+}$ in $M_{n,r}(\mathbb{R})$ be a semi-martingale. Then we
denote by
$\int_0^t A_{s-}\,dL_s B_{s-}$ the matrix $C_t$ in $M_{m,s}(\mathbb{R})$
which has ${ij}$th element $C_{ij,t}=\sum_{k=1}^n\sum_{l=1}^r\int
_0^t A_{ik,s-}B_{lj,s-}\,dL_{kl,s}$. Equivalently such an integral can be
understood in the sense of \cite{Metivier1982,Metivieretal1980} by
identifying it with the integral $\int_0^t \mathbf{A}_{s-}\,dL_s$ with
$\mathbf{A}_t$ being for each fixed $t$ the linear operator
$M_{n,r}(\mathbb{R})\to M_{m,s}(\mathbb{R}),  X\mapsto A_t X B_t$.
Analogous notation is used in the context of integrals with respect to
random measures.

Finally, integrals of the form $\int_A\int_B f(x,y) m(dx,dy)$ are
understood to be over the set $A$ in $x$ and over $B$ in $y$.

\subsection{L\'{e}vy bases}\label{sec:21}

To lay the foundations for the definition of vector-valued supOU
processes, we give now a summary of L\'{e}vy bases and the related
integration theory. In this context recall that a $d$-dimensional L\'
evy process can be understood as an $\mathbb{R}^d$-valued random
measure on
the real numbers. If $L=(L_t)_{t\in\mathbb{R}}$ is a $d$-dimensional
L\'evy
process, this measure is simply determined by $L((a,b])=L(b)-L(a)$ for
all $a,b\in\mathbb{R}, a<b$.

Define now $M_d^-:=\{X\in M_d(\mathbb{R})\dvtx \sigma(X)\subset(-\infty
,0)+i\mathbb{R}\}$ and $\mathscr{B}_b(M_d^-\times\mathbb
{R})$ to
be the bounded Borel sets of $M_d^-\times\mathbb{R}$. Note that
$M_d^-$ is
obviously a cone, but not a convex one (cf. \cite{Hershkowitz1998},
for instance). Moreover, we obviously have $\overline{M_d^-}=\{X\in
M_d(\mathbb{R})\dvtx \sigma(X)\subset(-\infty,0]+i\mathbb{R}\}$.

\begin{Definition}
A family $\Lambda=\{\Lambda(B):B\in\mathscr{B}_b(M_d^-\times
\mathbb{R})\}$ of $\mathbb{R}^d$-valued random variables is
called an
\emph{$\mathbb{R}^d$-valued L\'evy basis on $M_d^-\times\mathbb
{R}$} if:
\begin{enumerate}[(a)]
\item[(a)]  the distribution of $\Lambda(B)$ is infinitely divisible for all
$B\in\mathscr{B}_b(M_d^-\times\mathbb{R})$,
\item[(b)]  for any $n\in\mathbb{N}$ and pairwise disjoint sets $B_1,\ldots,
B_n\in\mathscr{B}_b(M_d^-\times\mathbb{R})$ the random
variables $\Lambda(B_1),\ldots,\Lambda(B_n)$ are independent and
\item[(c)]  for any pairwise disjoint sets $B_i\in\mathscr{B}_b
(M_d^-\times\mathbb{R})$, $i\in\mathbb{N}$, satisfying
$\bigcup_{n\in
\mathbb{N}} B_n\in\mathscr{B}_b(M_d^-\times\mathbb{R}
)$ the
series $\sum_{n=1}^\infty\Lambda(B_n)$ converges almost surely and
it holds that $\Lambda(\bigcup_{n\in\mathbb{N}} B_n
)=\sum
_{n=1}^\infty\Lambda(B_n)$.
\end{enumerate}
\end{Definition}

In the literature L\'evy bases are also often called {infinitely
divisible independently scattered random measures} (abbreviated
i.d.i.s.r.m.) instead.

In the following we will only consider L\'evy bases, which are
homogeneous (in time) and factorizable (into the effects of one
underlying infinitely divisible distribution and a probability
distribution on $M_d^-$); that is, their characteristic function is of
the form
%
\begin{equation}\label{eq:lbcf}
E(\exp(iu^*\Lambda(B)))=\exp(\varphi(u)\Pi
(B))
\end{equation}
for all $u\in\mathbb{R}^d$ and $B\in\mathscr{B}_b
(M_d^-(\mathbb{R}
)\times\mathbb{R})$. Here $\Pi=\pi\times\lambda$ is the product
of a probability measure $\pi$ on $M_d^-(\mathbb{R})$ and the Lebesgue
measure $\lambda$ on $\mathbb{R}$ and
%
\begin{equation}
\varphi(u)=iu^*\gamma-\frac{1}{2}u^*\Sigma u+\int_{\mathbb
{R}^d}
\bigl(e^{iu^*x}-1-iu^ *x1_{[0,1]}(\|x\|)\bigr)\nu(dx)
\end{equation}
is the cumulant transform of an infinitely divisible distribution on
$\mathbb{R}^d$ with L\'evy--Khintchine triplet $(\gamma,\Sigma,\nu)$,
that is, $\gamma\in\mathbb{R}^d$, $\Sigma\in\mathbb{S}_d^+$ and
$\nu$ is a
L\'evy measure---a~Borel measure on $\mathbb{R}^d$ with $\nu(\{0\}
)=0$ and
$\int_{\mathbb{R}^d}(\|x\|^2\wedge1)\nu(dx)<\infty$. The quadruple
$(\gamma,\Sigma,\nu,\pi)$ determines the distribution of the L\'evy
basis completely and is henceforth referred to as the ``generating
quadruple'' (cf. \cite{Fasen2005}).

The L\'evy process $L$ defined by
\[
L_t=\Lambda\bigl(M_d^-\times(0,t]\bigr)\quad\mbox{and}\quad L_{-t}=\Lambda
\bigl(M_d^-\times(-t,0)\bigr)\qquad\mbox{for } t\in\mathbb{R}^+,
\]
has characteristic triplet $(\gamma,\Sigma, \nu)$ and is called
``the underlying L\'evy process.''

For more information on $\mathbb{R}^d$-valued L\'evy bases see \cite
{Pedersen2003} and \cite{Rajputetal1989}.

A L\'{e}vy basis has a L\'{e}vy--It\^o decomposition.

\begin{Theorem}[\textup{(L\'evy--It\^o decomposition)}]\label{thm2.2}
Let $\Lambda$ be a homogeneous and factorisable $\mathbb{R}^d$-valued
L\'
evy basis on $M_d^-\times\mathbb{R}$ with generating quadruple
$(\gamma
,\Sigma,\nu,\pi)$. Then there exists a modification $\tilde\Lambda
$ of $\Lambda$ which is also a L\'evy basis with generating quadruple
$(\gamma,\Sigma,\nu,\pi)$ such that there exists an $\mathbb{R}
^d$-valued L\'evy basis $\tilde\Lambda^G$ on $M_d^-\times\mathbb
{R}$ with
generating quadruple $(0,\Sigma,0,\pi)$ and an independent Poisson
random measure $\mu$ on $(\mathbb{R}^d\times M_d^-\times\mathbb
{R},\mathscr
{B}(\mathbb{R}^d\times M_d^-\times\mathbb{R}))$ with intensity
measure $\nu
\times\pi\times\lambda$ which satisfy
\begin{eqnarray}\label{eqlevyito}
\tilde\Lambda(B)&=&\gamma(\pi\times\lambda)(B)+\tilde\Lambda
^G(B)\nonumber
\\
&&{}+\int_{\|x\|\leq1}\int_B x\bigl(\mu(dx,dA,ds)-ds\,\pi(dA)\nu
(dx)\bigr)
\\
&&{}+\int_{\|x\|> 1}\int_B x\mu(dx,dA,ds)\nonumber
\end{eqnarray}
for all $B\in\mathscr{B}_b(M_d^-\times\mathbb{R})$ and all $\omega
\in
\Omega$.

Provided $\int_{\|x\|\leq1}\|x\|\nu(dx)<\infty$, it holds that
\begin{eqnarray*}
\tilde\Lambda(B)=\gamma_0 (\pi\times\lambda)(B)+\tilde\Lambda
^G(B)+\int_{\mathbb{R}^d}\int_B x\mu(dx,dA,ds)
\end{eqnarray*}
for all $B\in\mathscr{B}_b(M_d^-\times\mathbb{R})$ with
%
\begin{equation}\label{eq:gamma0}
\gamma_0:=\gamma-\int_{\|x\|\leq1} x\nu(dx).
\end{equation}
Furthermore, the integral with respect to $\mu$ exists as a Lebesgue
integral for all $\omega\in\Omega$.
\end{Theorem}

Here an $\mathbb{R}^d$-valued L\'evy basis $\tilde\Lambda$ on
$M_d^-\times
\mathbb{R}$ is called a \emph{modification} of a L\'evy basis
$\Lambda$ if
$\tilde\Lambda(B)=\Lambda(B)$ a.s. for all $B\in\mathscr
{B}_b(M_d^-\times\mathbb{R})$. For the necessary background on the
integration with respect to Poisson random measures we refer to \cite
{Jacodetal2003}, Section 2.1, and \cite{Kallenberg2002}, Lemma 12.13.

\begin{pf*}{Proof of Theorem \ref{thm2.2}}
This follows immediately from \cite{Pedersen2003}, Theorem~4.5,
because the control measure $m$ is given by $m(B)=(\|\gamma\|
+\operatorname{tr}(\Sigma)+\int_{\mathbb{R}^d}(1\wedge\|x\|^2)\nu
(dx))(\pi
\times\lambda)(B)$ which is trivially continuous due to the presence
of the Lebesgue measure. The second part is an immediate consequence,
as no compensation for the small jumps is needed if $\int_{\|x\|\leq
1}\|x\|\nu(dx)<\infty$.
\end{pf*}

From now on we assume without loss of generality that all L\'evy bases
are such that they have the L\'evy--It\^o decomposition \eqref{eqlevyito}.

In the following we need to define integrals of deterministic functions
with respect to a L\'evy basis $\Lambda$. Following \cite
{Rajputetal1989}, for simple functions $f\dvtx M_d^ {-}\times\mathbb{R}\to
M_d(\mathbb{R})$,
\[
f(x)=\sum_{i=1}^m a_i 1_{B_i}(x)
\]
with $m\in\mathbb{N}$, $a_i\in M_d(\mathbb{R})$ and $B_i\in\mathscr
{B}_b
(M_d^-\times\mathbb{R})$, and for every $B\in\mathscr
{B}
(M_d^-\times\mathbb{R})$
we define the integral
\[
\int_{B} f(x)\Lambda(dx)=\sum_{i=1}^m a_i \Lambda(B\cap B_i).
\]
A $\mathscr{B}(M_d^-\times\mathbb{R})$-$\mathscr
{B}
(M_d(\mathbb{R}))$-measurable function $f\dvtx M_d^{-}\times\mathbb
{R}\to
M_d(\mathbb{R})$ is said to be \emph{$\Lambda$-integrable} if there exists
a sequence of simple functions $(f_n)_{n\in\mathbb{N}}$ such that
$f_n\to
f$ Lebesgue almost everywhere, and for all $B\in\mathscr{B}
(M_d^-\times\mathbb{R})$ the sequence $\int_B f_n(x)\Lambda(dx)$
converges in probability. For $\Lambda$-integrable $f$ we set $\int_B
f(x)\Lambda(dx)=\operatorname{plim}_{n\to\infty} \int
_Bf_n(x)\Lambda(dx)$. As
in \cite{Rajputetal1989} well definedness of the integral is ensured
by \cite{Urbaniketal1967}.

The following result is a straightforward generalization of \cite
{Rajputetal1989}, Propositions~2.6 and 2.7, to $\mathbb{R}^d$-valued L\'evy
bases and follows along the same lines.

\begin{Proposition}\label{th:exint}
Let $\Lambda$ be an $\mathbb{R}^d$-valued L\'{e}vy basis with
characteristic function of the form \eqref{eq:lbcf} and $f\dvtx M_d^
{-}\times\mathbb{R}\to M_d(\mathbb{R})$ a $\mathscr{B}
(M_d^-\times\mathbb{R}
)$-$\mathscr{B}(M_d(\mathbb{R}))$-meas\-ur\-able function.
Then $f$ is $\Lambda$-integra\-ble if and only if\
\begin{eqnarray}
&&\int_{M_d^-}\int_{\mathbb{R}}\biggl\|f(A,s)\gamma\nonumber
\\
&&{}\qquad\hspace*{13pt}+\int
_{\mathbb{R}^d}f(A,s)x\bigl(1_{[0,1]}(\|f(A,s)x\|)
\\
&&{}\hspace*{69pt}\hspace*{27pt}\hspace*{35pt}-1_{[0,1]}(\|x\|)\bigr)\nu(dx)\biggr\|\,ds\,\pi(dA)<\infty,\nonumber\vspace*{-6pt}
\end{eqnarray}
\begin{eqnarray}
\int_{M_d^-}\int_{\mathbb{R}} \|f(A,s)\Sigma f(A,s)^*\|\,ds\,\pi
(dA)&<&\infty
,\\
\int_{M_d^-}\int_{\mathbb{R}}\int_{\mathbb{R}^d}(1\wedge\|
f(A,s)x\|^2)\nu
(dx)\,ds\,\pi(dA)&<&\infty.
\end{eqnarray}

If $f$ is $\Lambda$-integrable, the distribution of $\int_{M_d^-}\int
_{\mathbb{R}^+}f(A,s)\Lambda(dA,ds)$ is infinitely divisible with
characteristic function
%
\begin{eqnarray}\label{eq:cfint}
&&E\biggl(\exp\biggl(iu^*\int_{M_d^-}\int_{\mathbb{R}^+}f(A,s)\Lambda
(dA,ds)\biggr)\biggr)\nonumber
\\[-8pt]\\[-8pt]
&&\qquad=\exp\biggl(\int_{M_d^-}\int_{\mathbb{R}^+}
\varphi(f(A,s)^*u)\,ds\,\pi(dA)\biggr)\nonumber
\end{eqnarray}
and characteristic triplet $(\gamma_{\mathrm{int}},\Sigma_{\mathrm
{int}},\nu_{\mathrm{int}})$ given by
\begin{eqnarray}
\gamma_{\mathrm{int}}&=&\int_{M_d^-}\int_{\mathbb{R}}
\biggl(f(A,s)\gamma\nonumber
\\
&&{}\hspace*{31pt}+\int_{\mathbb{R}^d}f(A,s)x\bigl(1_{[0,1]}(\|f(A,s)x\|
)
\\
&&{}\hspace*{130pt}-1_{[0,1]}(\|x\|)\bigr)\nu(dx)\biggr)\,ds\,\pi(dA),\nonumber
\\
\Sigma_{\mathrm{int}}&=&\int_{M_d^-}\int_{\mathbb{R}} f(A,s)\Sigma
f(A,s)^*\,ds\,\pi(dA),
\\
\qquad\qquad\nu_{\mathrm{int}}(B)&=&\int_{M_d^-}\int_{\mathbb{R}}\int_{\mathbb{R}
^d}1_B(f(A,s)x)\nu(dx)\,ds\,\pi(dA) \qquad \forall  B\in\mathscr{B}
(\mathbb{R}^d).
\end{eqnarray}
\end{Proposition}

When the underlying L\'evy process has finite variation we can do
$\omega$-wise Lebesgue integration; that is, the integral can be
obtained as a Lebesgue integral for each $ \omega\in\Omega$.

\begin{Proposition}\label{th:exlebsti}
Let $\Lambda$ be an $\mathbb{R}^d$-valued L\'{e}vy basis with
characteristic quadruple $(\gamma,0,\nu,\pi)$ satisfying $\int_{\|
x\|\leq1}\|x\|\nu(dx)<\infty$, and define $\gamma_0$ as in \eqref
{eq:gamma0},  that is, $\varphi(u)=iu^*\gamma_0+\int_{\mathbb
{R}^d}
(e^{iu^*x}-1)\nu(dx)$. Furthermore, let $f\dvtx M_d^ {-}\times
\mathbb{R}
\to M_d(\mathbb{R})$ be a $\mathscr{B}(M_d^-\times\mathbb
{R}
)$-$\mathscr{B}(M_d(\mathbb{R}))$-measurable function satisfying
\begin{eqnarray}\label{eq:finvarcond1}
\int_{M_d^-}\int_{\mathbb{R}}\|f(A,s)\|\,ds\,\pi(dA)&<&\infty,
\\\label{eq:finvarcond2}
\int_{M_d^-}\int_{\mathbb{R}}\int_{\mathbb{R}^d}\bigl(1\wedge\|
f(A,s)x\|\bigr)\nu
(dx)\,ds\,\pi(dA)&<&\infty.
\end{eqnarray}

Then
\begin{eqnarray}
\qquad\quad \int_{M_d^-}\int_{\mathbb{R}}f(A,s)\Lambda(dA,ds)\nonumber
&=&\int
_{M_d^-}\int
_{\mathbb{R}}f(A,s)\gamma_0 \,ds\,\pi(dA)
\\[-8pt]\\[-8pt]
&&{} +\int_{\mathbb{R}^d}\int
_{M_d^-}\int
_{\mathbb{R}}f(A,s) x\mu(dx,dA,ds),\nonumber
\end{eqnarray}
and the right-hand side is a Lebesgue integral for every $\omega\in
\Omega$ [conditions \eqref{eq:finvarcond1} and \eqref
{eq:finvarcond2} are also necessary for this].

Moreover, the distribution of $\int_{M_d^-}\int_{\mathbb
{R}}f(A,s)\Lambda
(dA,ds)$ is infinitely divisible with characteristic function
\begin{eqnarray*}
&&E\biggl(\exp\biggl(iu^*\int_{M_d^-}\int_{\mathbb{R}}f(A,s)\Lambda
(dA,ds)\biggr)\biggr)
\\
&&\qquad =e^{iu^*\gamma_{\mathrm{int},0}+\int_{\mathbb{R}
^d}(e^{iu^*x}-1)\nu_{\mathrm{int}}(dx)},\qquad u\in\mathbb{R}^d,
\end{eqnarray*}
where
\begin{eqnarray}
\quad\gamma_{\mathrm{int},0}&=&\int_{M_d^-}\int_{\mathbb{R}}f(A,s)\gamma
_0\,ds\,\pi(dA),
\\
\quad\qquad \nu_{\mathrm{int}}(B)&=&\int_{M_d^-}\int_{\mathbb{R}}\int_{\mathbb{R}
^d}1_B(f(A,s)x)\nu(dx)\,ds\,\pi(dA) \qquad \forall  B\in\mathscr
{B}(\mathbb{R}^d).
\end{eqnarray}
\end{Proposition}

\begin{pf}
Follows from the L\'evy--It\^o decomposition and the usual integration
theory with respect to Poisson random measures (see \cite{Kallenberg2002}, Lemma~12.13).
\end{pf}

\begin{Remark}
All results of this section remain valid when replacing $M_d^-$ with
$M_k(\mathbb{R})$, $k\in\mathbb{N}$, or any measurable subset of a finite-dimensional real vector space and when considering integration of
functions $f\dvtx M_k(\mathbb{R})\times\mathbb{R}\to M_{m,d}(\mathbb
{R})$. We decided to
state all our results with $M_d^-$ as this set will be used mainly in
the following and it reduces the notational burden.
\end{Remark}
%
\section{Multidimensional supOU processes}\label{sec:3}

In this section we introduce sup\-OU processes taking values in
$\mathbb{R}
^d$ with $d\in\mathbb{N}$ and analyze their properties. This extends
to a
multivariate setting the theory of univariate supOU processes as
introduced in \cite{barndorffnielsen01} and studied further, for
example, in \cite{FasenetCklu2007}.

Intuitively supOU processes are obtained by ``adding up'' independent
OU-type processes with different mean reversion coefficient.

\subsection{Definition and existence}\label{sec:31}
We define a $d$-dimensional supOU process as a process of the form
\eqref{eq:defsupou} below.

\begin{Theorem}\label{th:supouex}
Let $\Lambda$ be an $\mathbb{R}^d$-valued L\'evy basis on
$M_d^-\times\mathbb{R}
$ with generating quadruple $(\gamma, \Sigma,\nu,\pi)$
satisfying
\begin{equation}\label{c3}
\int_{\|x\|>1}\ln(\|x\|)\nu(dx)<\infty,
\end{equation}
and assume there exist measurable functions $\rho\dvtx M_d^-\to\mathbb{R}
^+\backslash\{0\}$ and $\kappa\dvtx  M_d^-\to[1,\infty)$ such that
\begin{eqnarray}\label{c1}
\|e^{As}\|\leq\kappa(A)e^{-\rho(A)s} \qquad \forall s\in\mathbb
{R}^+,
\pi\mbox{-almost surely},
\end{eqnarray}
 and
 \begin{eqnarray}\label{c2}
 \int_{M_d^-}\frac{\kappa(A)^2}{\rho(A)}\pi
(dA)<\infty.
\end{eqnarray}
Then the process $(X_t)_{t\in\mathbb{R}}$ given by
%
\begin{equation}\label{eq:defsupou}
X_t=\int_{M_d^-}\int_{-\infty}^te^{A(t-s)}\Lambda(dA,ds)
\end{equation}
is well defined for all $t\in\mathbb{R}$ and stationary. The distribution
of $X_t$ is infinitely divisible with characteristic triplet $(\gamma
_X,\Sigma_X,\nu_X)$ given by
\begin{eqnarray}
\gamma_X&=&\int_{M_d^-}\int_{\mathbb{R}^+}\biggl(e^{As}\gamma+\int
_{\mathbb{R}
^d}e^{As}x\bigl(1_{[0,1]}(\|e^{As}x\|)\nonumber
\\[-8pt]\\[-8pt]
&&{}\hspace*{127pt}-1_{[0,1]}(\|x\|)\bigr)\nu
(dx)\biggr)\,ds\,\pi(dA),\nonumber
\\
\Sigma_X&=&\int_{M_d^-}\int_{\mathbb{R}^+} e^{As}\Sigma
e^{A^*s}\,ds\,\pi
(dA),
\\
\qquad \nu_X(B)&=&\int_{M_d^-}\int_{\mathbb{R}^+}\int_{\mathbb
{R}^d}1_B(e^{As}x)\nu
(dx)\,ds\,\pi(dA)\nonumber
\\[-8pt]\\[-8pt]
&&\eqntext{\mbox{for all Borel sets } B\subseteq\mathbb{R}^d.}
\end{eqnarray}
\end{Theorem}
%

\begin{pf}
The stationarity is obvious once the well definedness is shown. Using
Proposition \ref{th:exint} it follows that necessary and sufficient
conditions for the integral to exist are given by
\begin{eqnarray}\label{eq:rjexcond1}
&&\int_{M_d^-}\int_{\mathbb{R}^+}\biggl\|e^{As}\gamma+\int_{\mathbb{R}
^d}e^{As}x\bigl(1_{[0,1]}(\|e^{As}x\|)\nonumber
\\[-8pt]\\[-8pt]
&&\qquad\hspace*{21pt}\hspace*{82pt}{} -1_{[0,1]}(\|x\|)\bigr)\nu
(dx)\biggr\|\,ds\,\pi(dA)<\infty,\nonumber
\\\label{eq:rjexcond2}
 &&\hspace*{41pt}\int_{M_d^-}\int_{\mathbb{R}^+} \|e^{As}\Sigma e^{A^*s}\|\,ds\,\pi
(dA)<\infty,
\\\label{eq:rjexcond3}
\qquad\quad &&\hspace*{-10pt}\int_{M_d^-}\int_{\mathbb{R}^+}\int_{\mathbb{R}^d}(1\wedge\|
e^{As}x\|^2)\nu
(dx)\,ds\,\pi(dA)<\infty.
\end{eqnarray}
First we show \eqref{eq:rjexcond3}
\begin{eqnarray*}
&&\int_{M_d^-}\int_{\mathbb{R}^+}\int_{\mathbb{R}^d}(1\wedge\|
e^{As}x\|^2)\nu
(dx)\,ds\,\pi(dA)
\\
&&\qquad\leq\int_{M_d^-}\int_{\mathbb
{R}^+}\int
_{\mathbb{R}^d}\bigl(1\wedge\kappa(A)^2e^{-2\rho(A)s}\|x\|^2\bigr)\nu
(dx)\,ds\,\pi
(dA)
\\
&&\qquad=\int_{M_d^-}\int_{\|x\|>1/\kappa(A)}\frac{\ln
(\kappa(A)\|x\|)+1/2}{\rho(A)}\nu(dx)\pi(dA)
\\
&&\quad\quad\quad{}+\int_{M_d^-}\int_{\|x\|\leq1/\kappa(A)}\frac
{\kappa(A)^2\|x\|^2}{2\rho(A)}\nu(dx)\pi(dA).\vspace*{2pt}
\end{eqnarray*}
The finiteness of the first integral follows from \eqref{c3}, \eqref
{c2}, $\kappa(A)\geq1$ and $\nu$ being a L\'evy measure, which imply\vspace*{2pt}
\begin{eqnarray*}
&&\int_{M_d^-}\int_{\|x\|>1/\kappa(A)}\frac{\ln(\kappa(A)\|x\|
)+1/2}{\rho(A)}\nu(dx)\pi(dA)
\\
&&\qquad\leq\int_{M_d^-}\int_{\|x\|>1}\frac{\ln(\kappa(A))+\ln
(\|x\|)+1/2}{\rho(A)}\nu(dx)\pi(dA)
\\
&&{}\qquad\quad+\int_{M_d^-}\int_{\|x\|\leq1}\frac{3\kappa(A)^2\|
x\|^2}{2\rho(A)}\nu(dx)\pi(dA)
\\
&&\qquad=\int_{M_d^-}\frac{\ln(\kappa(A))}{\rho(A)}\pi(dA)\int
_{\|x\|>1}\nu(dx)
\\
&&{}\quad\qquad+\int_{M_d^-}\frac{3\kappa(A)^2}{2\rho(A)}\pi
(dA)\int_{\|x\|\leq1}\|x\|^2\nu(dx)
\\
&&{}\qquad\quad+\int_{M_d^-}\frac{1}{\rho(A)}\pi(dA)\int_{\|x\|
>1}\bigl(\ln(\|x\|)+1/2\bigr)\nu(dx)<\infty.\vspace*{2pt}
\end{eqnarray*}

Likewise the finiteness of the second integral is implied by \eqref
{c2}, $\kappa(A)\geq1$ and $\int_{\|x\|\leq1}\|x\|^2\nu(dx)<\infty
$, as $\nu$ is a L\'evy measure.

Next \eqref{eq:rjexcond2} follows from \eqref{c2} and
\begin{eqnarray*}
\int_{M_d^-}\int_{\mathbb{R}^+} \|e^{As}\Sigma e^{A^*s}\|\,ds\,\pi
(dA)&\leq&\|
\Sigma\|\int_{M_d^-}\int_{\mathbb{R}^+} \kappa(A)^2e^{-2\rho
(A)s}\,ds\,\pi(dA)
\\
&=&\|\Sigma\|\int_{M_d^-}\frac{\kappa(A)^2}{2\rho(A)}\pi(dA).
\end{eqnarray*}

Turning to \eqref{eq:rjexcond1} we have from \eqref{c2} that
\begin{eqnarray*}
\int_{M_d^-}\int_{\mathbb{R}^+}\|e^{As}\gamma\|\,ds\,\pi
(dA)&\leq&\|\gamma\|\int_{M_d^-}\int_{\mathbb{R}^+}\kappa
(A)e^{-\rho
(A)s}\,ds\,\pi(dA)
\\
&=&\|\gamma\|\int_{M_d^-}\frac{\kappa(A)}{\rho
(A)}\pi(dA)<\infty.
\end{eqnarray*}
Moreover,
\begin{eqnarray*}
&&\int_{M_d^-}\int_{\mathbb{R}^+}\biggl\|\int_{\mathbb
{R}^d}e^{As}x
\bigl(1_{[-1,1]}(\|e^{As}x\|)-1_{[-1,1]}(\|x\|)\bigr)\nu(dx)\biggr\|
\,ds\,\pi(dA)
\\
&&\qquad\leq\int_{M_d^-}\int_{\mathbb{R}^+}\int_{\|x\|\leq1, \|
e^{As}x\|
\geq1}\|e^{As}x\|\nu(dx)\,ds\,\pi(dA)
\\
&&{}\qquad\quad+\int
_{M_d^-}\int_{\mathbb{R}^+}\int_{\|x\|\geq1, \|e^{As}x\|\leq
1}\|
e^{As}x\|\nu(dx)\,ds\,\pi(dA)
\\
&&\qquad\leq\int_{M_d^-}\int_{\mathbb{R}^+}\int_{\|x\|\leq1, \|
e^{As}x\|
\geq1}\|e^{As}x\|^2\nu(dx)\,ds\,\pi(dA)
\\
&&{}\qquad\quad+\int_{M_d^-}\int_{\mathbb{R}^+}\int_{\|x\|\in
(1,e^{\rho(A)s/2})} \|x\|\kappa(A)e^{-\rho(A) s}\nu(dx)\,ds\,\pi(dA)
\\
&&{}\qquad\quad+\int_{M_d^-}\int_{\mathbb{R}^+}\int_{\|x\|\geq
e^{\rho(A)s/2}} \nu(dx)\,ds\,\pi(dA)
\\
&&\qquad\leq\int_{\|x\|\leq1}\|x\|^2\nu(dx)\int_{M_d^-}\frac
{\kappa(A)^2}{2\rho(A)}\pi(dA)+ \int_{M_d^-}\frac{2\kappa
(A)}{\rho(A)}\pi(dA)\int_{\|x\|>1}\nu(dx)
\\
&&{}\qquad\quad+\int_{M_d^-}\frac{2}{\rho(A)}\pi(dA)\int_{\|x\|
>1}\ln(\|x\|)\nu(dx)<\infty
\end{eqnarray*}
with the finiteness following from \eqref{c3}, \eqref{c2} and $\nu$
being a L\'evy measure.

That the distribution of $X_t$ is infinitely divisible and has the
stated characteristic triplet follows now immediately from Proposition
\ref{th:exint}.
\end{pf}

\begin{Remark}
(i) The necessary and sufficient conditions for the existence of the
multivariate supOU process $X$ are \eqref{eq:rjexcond1}--\eqref{eq:rjexcond3}. However, as they are obviously
very intricate to check in concrete situations, it seems to be
appropriate to replace them by the sufficient
conditions \eqref{c3}--\eqref{c2}. One particular advantage of these conditions
is that they involve only integrals with respect to either $\nu$ or
$\pi$, but not with respect to both.\vspace*{-6pt}
\begin{longlist}[(iii)]
\item[(ii)] Note also that for $d=1$ the conditions above become the necessary
and sufficient conditions of \cite{FasenetCklu2007}, as we can then
take $\kappa(A)=1$ and $\rho(A)=-A$.

\item[(iii)] By looking at the Jordan decomposition one can see that pointwise
there is for any $A\in M_d^-$ a constant $\kappa\in[1,\infty)$ and a
$\rho\in(0,-\max(\Re(\sigma(A)))]$ such that $\|e^{As}\|\leq
\kappa e^{-\rho s}$ for all $s\in\mathbb{R}^+$. If $A$ is diagonalizable,
it is possible to choose $\rho(A)=-\max(\Re(\sigma(A)))$ and
$\kappa(A)=\|U\|\|U^{-1}\|$ with $U\in GL_d(\mathbb{C})$ being such that
$UAU^{-1}$ is diagonal. So \eqref{c1} essentially demands that this
choice has to be done measurably in $A$ (see especially Example \ref
{ex5} for a concrete example).
\end{longlist}
\end{Remark}

It seems important, especially for applications, to understand how
close our sufficient existence conditions are to necessity and to give
also necessary conditions easier checkable than \eqref{eq:rjexcond1}--\eqref{eq:rjexcond3}. To this end we need the
concept of modulus of injectivity (see, e.g., \cite{Pietsch1980}, Section B.3). For a $Z\in M_d(\mathbb{R})$ the modulus of
injectivity is
defined as $j(Z)=\min_{\|x\|=1}\|Zx\|$. It is immediate to see that
$0\leq j(Z)\leq\|Z\|$ and $j(Z)=\|Z^{-1}\|^{-1}$ for $Z\in
GL_d(\mathbb{R})$.
\begin{Proposition}\label{th:vexcond}
Let $\Lambda$ be an $\mathbb{R}^d$-valued L\'evy basis on
$M_d^-\times\mathbb{R}
$ with generating quadruple $(\gamma, \Sigma,\nu,\pi)$
and assume there exist measurable functions $\tau:M_d^-\to\mathbb{R}
^+\backslash\{0\}$ and $\vartheta\dvtx M_d^-\to(0,1]$ such that
\begin{eqnarray}\label{nc1}
j(e^{As})&\geq\vartheta(A)e^{-\tau(A)s}\qquad  \forall s\in\mathbb{R}
^+,\pi\mbox{-almost surely}.
\end{eqnarray}
Then necessary conditions for the integral \eqref{eq:defsupou} to
exist are
\begin{eqnarray}\label{nc2a}
 &&\qquad\quad\int_{\vartheta(A)\geq\varepsilon}\frac{1}{\tau(A)}\pi(dA)<\infty
\nonumber
\\[-8pt]\\[-8pt]
&&\qquad\qquad\ \mbox{for any }\varepsilon\in(0,1] \mbox{ with }\nu\bigl(\{\|
x\|> 1/\varepsilon\}\bigr)>0
 \mbox{ and } \pi(\{\vartheta
(A)\geq\varepsilon\})>0,\nonumber
\\\label{nc2b}
&&\hspace*{32pt}\int_{M_d^-}\frac{\vartheta(A)^2}{\tau(A)}\pi(dA)<\infty\qquad \mbox{provided }j(\Sigma)>0\mbox{ or }\nu(\{\|x\|\leq1\})>0
\end{eqnarray}
and
\begin{eqnarray}\label{nc3}
 \int_{\|x\|> 1}\ln(\|x\|)\nu(dx)<\infty. 
\end{eqnarray}
\end{Proposition}

This shows, in particular, that the logarithmic moment condition on
$\nu$ is both necessary and sufficient and that \eqref{c1}, \eqref
{c2} together with \eqref{nc1}--\eqref{nc2b} form a
set of sufficient and necessary conditions which are as close as one
can probably hope for conditions reasonably easy to work with.

\begin{pf*}{Proof of Proposition \ref{th:vexcond}}
In the case $j(\Sigma)>0$ condition \eqref{nc2b} follows immediately
from \eqref{eq:rjexcond2}.

Equation \eqref{eq:rjexcond3} implies that a necessary condition is\vspace*{-1pt}
\begin{eqnarray*}
\infty&>& \int_{M_d^-}\int_{\mathbb{R}^+}\int_{\mathbb
{R}^d}\bigl(1\wedge\vartheta
(A)^2e^{-2\tau(A)s}\|x\|^2\bigr)\nu(dx)\,ds\,\pi(dA)
\\[-2pt]
&=&\int_{M_d^-}\int
_{\|x\|>1/\vartheta(A)}\frac{\ln(\vartheta(A)\|x\|)+1/2}{\tau
(A)}\nu(dx)\pi(dA)
\\[-2pt]
&&{}+\int_{M_d^-}\int_{\|x\|\leq1/\vartheta(A)}\frac
{\vartheta(A)^2\|x\|^2}{2\tau(A)}\nu(dx)\pi(dA).
\end{eqnarray*}

\noindent
Since $\vartheta(A)\leq1$, the second summand implies \eqref{nc2b}
given $\nu(\{\|x\|\leq1\})>0$.

Choose now $\varepsilon>0$ such that $\pi(\{\vartheta(A)\geq\varepsilon\}
)>0$. Then the first integral is bigger than
\begin{eqnarray*}
\int_{\vartheta(A)\geq\varepsilon}\int_{\|x\|>1/\varepsilon}\frac{\ln
(\vartheta(A)\|x\|)+1/2}{\tau(A)}\nu(dx)\pi(dA).
\end{eqnarray*}
This gives the necessity of \eqref{nc2a}. Moreover, the above integral
is again greater than
\
\begin{eqnarray*}
\int_{\vartheta(A)\geq\varepsilon}\int_{\|x\|>1/\varepsilon}\frac{\ln
(\varepsilon)+\ln(\|x\|)+1/2}{\tau(A)}\nu(dx)\pi(dA).
\end{eqnarray*}
Since $\int_{\vartheta(A)\geq\varepsilon}\int_{\|x\|>1/\varepsilon
}\frac{1}{\tau(A)}\nu(dx)\pi(dA)$ is necessarily finite, we have
that
\[
\int_{\vartheta(A)\geq\varepsilon}\int_{\|x\|>1/\varepsilon}\frac{\ln
(\varepsilon)}{\tau(A)}\nu(dx)\pi(dA)>-\infty.
\]
Hence, $\int_{\vartheta(A)\geq\varepsilon}\int_{\|x\|>1/\varepsilon
}\frac{\ln(\|x\|)}{\tau(A)}\nu(dx)\pi(dA)<\infty$. Since\vspace*{1pt} $\int
_{\vartheta(A)\geq\varepsilon}\frac{1}{\tau(A)}\pi (dA)>0$ by
construction, this implies the necessity of \eqref{nc3}.
\end{pf*}

\begin{Remark}
(i) For $d=1$ we have again recovered the necessary and sufficient
conditions of \cite{FasenetCklu2007}, as we can then take $\vartheta
(A)=1$ and $\tau(A)=-A$.\vspace*{-6pt}
\begin{longlist}[(iii)]
\item[(ii)] Pointwise there is again for any $A\in M_d^-$ a constant
$\vartheta\in(0,1]$ and a $\tau\in[-\min(\Re(\sigma(A))),\infty
)$ such that $j(e^{As})\geq\vartheta e^{-\tau s}$ for all $s\in
\mathbb{R}
^+$. For $A$ diagonalizable $\tau(A)=-\min(\Re(\sigma(A)))$ and
$\vartheta(A)=j(U)j(U^{-1})$ can be chosen if $UAU^{-1}$ is diagonal.

\item[(iii)] If $\nu$ has unbounded support, $\int_{\vartheta(A)\geq
\varepsilon}\frac{1}{\tau(A)}\pi(dA)<\infty$ is necessary for all
$\varepsilon>0$. If $\pi(\{\vartheta(A)\geq\tilde\varepsilon\})=1$ for
some $\tilde\varepsilon>0$, then $\int_{M_d^{-}}\frac{1}{\tau(A)}\pi
(dA)<\infty$ becomes a necessary condition, provided $\nu(\{\|x\|
>1/\tilde\varepsilon\})>0$.
\end{longlist}
\end{Remark}

In some applications like stochastic volatility modeling, for instance,
one is particularly interested in the case where the underlying L\'evy
process is of finite variation and the supOU process is defined via
$\omega$-wise integration. The following result is proved using
Proposition \ref{th:exlebsti} together with variations of the
arguments of the proofs of Theorem \ref{th:supouex} and Proposition
\ref{th:vexcond}.

\begin{Proposition}\label{th:supoufv}
(i) Let $\Lambda$ be an $\mathbb{R}^d$-valued L\'evy basis on
$M_d^-\times
\mathbb{R}$ with generating quadruple $(\gamma, 0,\nu,\pi)$ satisfying
\begin{equation}\label{fvc3}
\int_{\|x\|>1}\ln(\|x\|)\nu(dx)<\infty \quad\mbox{and}\quad   \int_{\|
x\|\leq1}\|x\|\nu(dx)<\infty
\end{equation}
and assume there exist measurable functions $\rho\dvtx M_d^-\to\mathbb{R}
^+\backslash\{0\}$ and $\kappa\dvtx M_d^-\to[1,\infty)$ such that
\begin{eqnarray}\label{fvc1}
\|e^{As}\|\leq\kappa(A)e^{-\rho(A)s} \qquad \forall s\in\mathbb
{R}^+,
\pi\mbox{-almost surely},
\end{eqnarray}
 and
\begin{eqnarray}\label{fvc2}
\int_{M_d^-}\frac{\kappa(A)}{\rho(A)}\pi
(dA)<\infty.
\end{eqnarray}
Then the process $(X_t)_{t\in\mathbb{R}}$ given by
\begin{eqnarray}\label{eq:lebint}
\qquad\hspace*{8pt} X_t&=&\int_{M_d^-}\int_{-\infty}^te^{A(t-s)}\Lambda(dA,ds)\nonumber
\\[-8pt]\\[-8pt]
&=&\int_{M_d^-}\int_{-\infty}^te^{A(t-s)}\gamma_0 \,ds\,\pi
(dA)+\int_{\mathbb{R}^d}\int_{M_d^-}\int_{-\infty}^te^{A(t-s)}
x\mu
(dx,dA,ds)\nonumber
\end{eqnarray}
is well defined as a Lebesgue integral for all $t\in\mathbb{R}$ and
$\omega
\in\Omega$ and $X$ is stationary.

\textup{(ii)} If there exist measurable functions $\tau\dvtx M_d^-\to\mathbb{R}
^+\backslash\{0\}$ and $\vartheta: M_d^-\to(0,1]$ such that
$
j(e^{As})\geq\vartheta(A)e^{-\tau(A)s}$ $ \forall s\in\mathbb{R}
^+,\pi\mbox{-almost surely},
$
then necessary conditions for the integral \eqref{eq:lebint} to exist
as a Lebesgue integral are:
\begin{eqnarray}\label{fvnc2a}
\qquad\nonumber\int_{\vartheta(A)\geq\varepsilon}\frac{1}{\tau(A)}\pi
(dA)&<& \infty
\\[-8pt]\\[-8pt]
\eqntext{\mbox{for any }\varepsilon\in(0,1]
\mbox{ such that }\nu(\{\|x\|> 1/\varepsilon\})>0 \mbox{ and } \pi\bigl(\{
\vartheta(A)\geq\varepsilon\}\bigr)>0,}
\\\label{fvnc2b}
\int_{M_d^-}\frac{\vartheta(A)}{\tau(A)}\pi(dA)&<&\infty
\qquad\mbox{provided }\gamma_0\neq 0\mbox{ or }\nu(\{\|x\|\leq1\})>0,
\\\label{fvnc3}
\int_{\|x\|> 1}\ln(\|x\|)\nu(dx)&<&\infty\quad\mbox{and}\quad \int_{\|x\|
\leq1}\|x\|\nu(dx)<\infty.
\end{eqnarray}
\end{Proposition}

\begin{Remark}
If \eqref{c2} is satisfied for a L\'evy basis, then \eqref{fvc2} is
also satisfied.

\end{Remark}

We shall not develop the general case further, but consider two special
cases which appear to be sufficient for most purposes.
We define $M_d^{N-}:=\{A\in M_d(\mathbb{R}): A \mbox{ is normal and }
\sigma(A)\subset(-\infty,0)+i\mathbb{R}\}.$

\begin{Proposition}\label{th:exprop}
\textup{(i)} Assume that $\pi(M_d^{N-})=1$, then \eqref{c1} or \eqref{fvc1}
are satisfied with $\kappa(A)=1$ and $\rho(A)=-\max(\Re(\sigma
(A)))$. Moreover, \eqref{c2} or \eqref{fvc2} are implied by
\begin{eqnarray}\label{eq:condexnormal}
-\int_{M_d^{N-}}\frac{1}{\max(\Re(\sigma(A)))}\pi(dA)&<\infty.
\end{eqnarray}

Likewise, \eqref{nc1} is satisfied with $\vartheta(A)=1$ and $\tau
(A)=-\min(\Re(\sigma(A)))$. The necessary conditions \eqref{nc2a},
\eqref{nc2b}, \eqref{fvnc2a} and \eqref{fvnc2b} all become
\[
-\int_{M_d^{N-}}\frac{1}{\min(\Re(\sigma(A)))}\pi(dA)<\infty
\]
[assuming $\nu(\|x\|>1)>0$ for \eqref{nc2a}, \eqref{fvnc2a}].

\textup{(ii)} Assume that there are a $K\in\mathbb{N}$ and diagonalizable
$A_1,\ldots,A_K\in M_d^-{(\mathbb{R})}$ such that $\pi(\{\lambda
A_i\dvtx
i=1,\ldots,K; \lambda\in\mathbb{R}^+\backslash\{0\}\})=1$. Then
\eqref
{c1} or \eqref{fvc1} are satisfied with $\kappa(A)=C$ for some $C\in
[1,\infty)$ and $\rho(A)=-\max(\Re(\sigma(A)))$.
Moreover, \eqref{c2} or \eqref{fvc2} are implied by
%
\begin{equation}\label{eq:condexray}
-\int_{M_d^{-}}\frac{1}{\max(\Re(\sigma(A)))}\pi(dA)<\infty.
\end{equation}
Likewise, \eqref{nc1} is satisfied with $\tau(A)=-\min(\Re(\sigma
(A)))$ and $\vartheta(A)=c$ for some $c\in(0,1]$ and the necessary
conditions \eqref{nc2a}, \eqref{nc2b}, \eqref{fvnc2a} and \eqref
{fvnc2b} all become
$
-\int_{M_d^{-}}\frac{1}{\min(\Re(\sigma(A)))}\pi(dA)<\infty
$ [assuming $\nu(\|x\|>1/c)>0$ for \eqref{nc2a}, \eqref{fvnc2a}].
\end{Proposition}

In dimension one these are again the well-known necessary and
sufficient conditions. Observe also that the eigenvalues are continuous
(and hence measurable) in $A$ because they are the zeros of the
characteristic polynomial.

\begin{pf*}{Proof of Proposition \ref{th:exprop}}
Part (i) follows immediately from the fact that all normal matrices are
unitarily diagonalizable.

Likewise, (ii) is a consequence of the above mentioned pointwise bound
and the fact that this can be turned into a global one because for
fixed $i=1,\ldots, N$ the matrices $\{\lambda A_i\}_{\lambda\in
\mathbb{R}
^+\backslash\{0\}}$ are all diagonalized by the same invertible matrices.
\end{pf*}

In (i) the mean reversion parameter $A$ of the superimposed OU-type
processes is restricted to normal matrices and in (ii) to finitely many
rays $\{\lambda A_i\}_{\lambda\in\mathbb{R}^+\backslash\{0\}}$.

\begin{Remark}
\textup{(i)} Typically one will, in general, not consider normal matrices for
$A$ as in \textup{(i)}, but only negative definite ones, since this allows one
to use well-known distributions on the positive definite matrices (see,
e.g., \cite{Guptaetal2000}) for $\pi$. In the case \textup{(ii)} possible $\pi
$ can be obtained by using arbitrary distributions on $\mathbb{R}^+$ along
the rays and positive weights summing to one for the different rays.

\textup{(ii)} Intuitively \eqref{eq:condexnormal} and \eqref{eq:condexray} and
their necessary counterparts mean that $\pi$ must not put too much
mass on elements with very slow exponential decay rates.
\end{Remark}

\subsection{Finiteness of moments and second-order structure}

Before we look at the second-order structure, we give conditions
ensuring the finiteness of moments.
\begin{Theorem}\label{th:exmom}
Let $X$ be a stationary $d$-dimensional supOU process driven by a L\'
evy basis $\Lambda$ satisfying the conditions of Theorem \ref{th:supouex}.\vspace*{-6pt}
\begin{longlist}[(iii)]
\item[(i)] If
%
\begin{equation}\label{eq:excondmoms2}
\int_{\|x\|>1}\|x\|^r\nu(dx)<\infty
\end{equation}
for $r\in(0,2]$, then $X$ has a finite $r$th moment, that is, $E(\|
X_t\|^r)<\infty$.

\item[(ii)] If $r\in(2,\infty)$ and
%
\begin{equation}\label{eq:excondmom}
\int_{\|x\|>1}\|x\|^r\nu(dx)<\infty,\qquad\int_{M_d^-}\frac{\kappa
(A)^r}{\rho(A)}\pi(dA)<\infty,
\end{equation}
then $X$ has a finite $r$th moment, that is, $E(\|X_t\|^r)<\infty$.

\item[(iii)] Necessary conditions for $X$ to have a finite $r$th moment are
%
\begin{equation}
\int_{\|x\|>1}\|x\|^r\nu(dx)<\infty
\end{equation}
in general and
%
\begin{equation}
\int_{\vartheta(A)\geq\varepsilon}\frac{\vartheta(A)^r}{\tau(A)}\pi
(dA)<\infty
\end{equation}
for any $\varepsilon$ such that $\nu(\{\|x\|> 1/\varepsilon\})>0 \mbox{
and } \pi(\{\vartheta(A)\geq\varepsilon\})>0$.
\end{longlist}
\end{Theorem}

In connection with the above results observe that the underlying L\'evy
process~$L$ has an $r$th moment, that is, $E(\|L_1\|^r)<\infty$, for
$r\in\mathbb{R}^+$ if and only if
\[
\int_{\|x\|>1}\|x\|^r \nu_L(dx)
\mbox{ is finite}.
\]

\begin{pf*}{Proof of Theorem \ref{th:exmom}}
Using \cite{Sato1999}, Corollary 25.8, we have to show $\int_{\|x\|>
1}\|x\|^r\nu_X(dx)<\infty$ to establish (i) and (ii).
Now,
\begin{eqnarray*}
&&\int_{\|x\|> 1}\|x\|^r\nu_X(dx)
\\
&&\qquad = \int_{M_d^-}\int_0^\infty\int
_{\mathbb{R}^d}\|e^{As}x\|^r1_{(1,\infty)}(\|e^{As}x\|)\nu(dx)\,ds\,\pi
(dA)
\\
&&\quad\quad\leq\int_{M_d^-}\int_0^\infty\int_{\mathbb
{R}^d}\kappa
(A)^re^{-r\rho(A)s}\|x\|^r1_{(1,\infty)}\bigl(\kappa(A)e^{-\rho(A)s}\|x\|
\bigr)\nu(dx)\,ds\,\pi(dA)\
\\
&&\quad\quad=\int_{M_d^-}\int_{\|x\|>1/\kappa(A)}\int_0^{
{\ln(\kappa(A)\|x\|)/\rho(A)}}\kappa(A)^re^{-r\rho(A)s}\|x\|
^r\,ds\,\nu(dx)\pi(dA)\
\\
&&\quad\quad=\int_{M_d^-}\int_{\|x\|>1/\kappa(A)}\frac{\kappa
(A)^r\|x\|^r}{r\rho(A)}\biggl(1-\frac{1}{\kappa(A)^r\|x\|^r}
\biggr)\nu(dx)\pi(dA)
\\
&&\quad\quad=\int_{M_d^-}\int_{\|x\|>1/\kappa(A)}\frac{\kappa
(A)^r\|x\|^r-1}{r\rho(A)}\nu(dx)\pi(dA).
\end{eqnarray*}
That $\int_{M_d^-}\int_{\|x\|>1/\kappa(A)}\frac{1}{r\rho(A)}\nu
(dx)\pi(dA)<\infty$ has already been shown in the proof of Theorem
\ref{th:supouex}.

Moreover, we obtain
\begin{eqnarray*}
&&\int_{M_d^-}\int_{\|x\|> 1/\kappa(A)}\frac{\kappa(A)^r\|x\|
^r}{\rho(A)}\nu(dx)\pi(dA)
\\
&&\qquad\leq\int_{M_d^-}\int_{\|x\|>
1}\frac{\kappa(A)^r\|x\|^r}{\rho(A)}\nu(dx)\pi(dA)
\\
&&{}\qquad\quad+\int
_{M_d^-}\int_{\|x\|\leq1}\frac{\kappa(A)^{r\vee2}\|x\|^{r\vee
2}}{\rho(A)}\nu(dx)\pi(dA).
\end{eqnarray*}
Hence, (i) and (ii) follow, since $\nu$ is a L\'evy measure, using
also \eqref{c2} for (i).

Regarding the proof of (iii), analogous arguments give that
\[
\int_{M_d^-}\int_{\|x\|>1/\vartheta(A)}\frac{\vartheta(A)^r\|x\|
^r-1}{r\tau(A)}\nu(dx)\pi(dA)<\infty
\]
is a necessary condition. An inspection of the proof of Proposition
\ref{th:vexcond} shows that
\[
\int_{M_d^-}\int_{\|x\|>1/\vartheta(A)}\bigl(1/\tau(A)\bigr)\nu(dx)\pi
(dA)<\infty
\]
is already necessary for $X$ to exist. Hence,
\begin{eqnarray*}
\infty&>& \int_{M_d^-}\int_{\|x\|>1/\vartheta(A)}\frac{\vartheta
(A)^r\|x\|^r}{r\tau(A)}\nu(dx)\pi(dA)
\\
&\geq&\int_{\vartheta
(A)\geq\varepsilon}\int_{\|x\|>1/\varepsilon}\frac{\vartheta(A)^r\|x\|
^r}{r\tau(A)}\nu(dx)\pi(dA)
\end{eqnarray*}
is necessary for a finite $r$th moment of $X$. This implies (iii),
since\break $\int_{\|x\|>1}\|x\|^r\nu(dx)$ is finite if and only if $\int
_{\|x\|>c}\|x\|^r\nu(dx)<\infty$ for arbitrary $c>0$.
\end{pf*}

\begin{Remark}
In the set-up of Proposition \ref{th:exprop}(i) [and analogously in (ii)]
\begin{eqnarray}\label{eq:condmoms}
\quad -\int_{M_d^{-}}\frac{1}{\max(\Re(\sigma(A)))}\pi(dA)<\infty
,\qquad\int_{\|x\|> 1}\|x\|^r\nu(dx)<\infty
\end{eqnarray}
imply \eqref{eq:condexnormal} and \eqref{eq:excondmom} [resp. \eqref
{eq:excondmoms2}].

Likewise, $
-\int_{M_d^{-}}\frac{1}{\min(\Re(\sigma(A)))}\pi(dA)<\infty$,
provided $\nu(\|x\|>1)>0$ or $\nu(\|x\|>1/c)>0$, respectively, and
$\int_{\|x\|> 1}\|x\|^r\nu(dx)<\infty$ become necessary conditions
for $X$ to exist and to have a finite $r$th moment.
\end{Remark}

In applications these results have important implications for modeling.
If one wants to have certain moments finite and certain moments
infinite, for example, because this is what observed data strongly
suggests, one has to use a driving L\'evy basis $\Lambda$ having
exactly the same moments finite.

Moreover, knowledge of the moments allows for the estimation of the
model based on empirical observations by using the general method of
moments (GMM) estimation procedure, for example. Often the inference on
parameters is based on the second-order moment structure as, for
instance, in the estimation of multivariate OU-type stochastic
volatility models in \cite{PigorschetStelzer2006} or the one of
univariate supOU type models in \cite{Barndorffetal2003}. To provide
the foundations for such work and due to the importance in applications
of understanding the temporal dependence structure, we now calculate
the second-order moment structure.

\begin{Theorem}\label{th:secord}
Let $X$ be a stationary $d$-dimensional supOU process driven by a L\'
evy basis $\Lambda$ satisfying the conditions of Theorem \ref
{th:supouex} and assume additionally that $\int_{\mathbb{R}^d}\|x\|
^2\nu
(dx)<\infty$. Then $E(\|X_0\|^2)<\infty$ and we have
\begin{eqnarray}\label{eq:exp}
E(X_0)&=&-\int_{M_d^-}A^{-1}\biggl(\gamma+\int_{|x|>1}x\nu(dx)
\biggr)\pi(dA),
\\\label{eq:var}
\operatorname{var}(X_0)&=&-\int_{M_d^-}(\mathscr{A}(A))^{-1}
\biggl(\Sigma+
\biggl(\int_{\mathbb{R}^d}xx^*\nu(dx)\biggr)\biggr)\pi(dA),
\\\label{eq:covar}
\operatorname{cov}(X_h,X_0)&=&-\int_{M_d^-}e^{Ah}(\mathscr
{A}(A))^{-1}\biggl(\Sigma
+\int_{\mathbb{R}^d}xx^*\nu(dx)\biggr)\pi(dA)\nonumber
\\[-8pt]\\[-8pt]
&&\eqntext{\mbox{ for } h\in
\mathbb{R}
^+}
\end{eqnarray}
with $\mathscr{A}(A)\dvtx M_d(\mathbb{R})\to M_d(\mathbb{R}), X\mapsto AX+XA^*$.

Moreover, it holds that
%
\begin{equation}\label{eq:acovXzero}
\lim_{h\to\infty}\operatorname{cov}(X_h,X_0)=0.
\end{equation}
\end{Theorem}

\begin{pf}
The finiteness of the second moments follows from Theorem \ref
{th:exmom}. Using the formulae of Theorem \ref{th:supouex} and \cite
{Sato1999}, Example 25.12, we obtain
\begin{eqnarray*}
E(X_0)=\gamma_X+\int_{\|x\|> 1}x\nu_X(dx)=\int_{M_d^-}\int
_{\mathbb{R}
^+}e^{As}\biggl(\gamma+\int_{\|x\|> 1}x\nu(dx)\biggr)\,ds\,\pi(dA).
\end{eqnarray*}
Noting that $\frac{d}{ds}A^{-1}e^{As}=e^{As}$, integrating over $s$
gives \eqref{eq:exp}.

Likewise we get
\begin{eqnarray*}
\operatorname{var}(X_0)&=&\Sigma_X+\int_{\mathbb{R}^d}xx^*\nu
_X(dx)
\\
&=&\int_{M_d^-}\int
_{\mathbb{R}^+}e^{As}\biggl(\Sigma+\int_{\mathbb{R}^d}xx^*\nu
(dx)\biggr)
e^{A^*s}\,ds\,\pi(dA)
\end{eqnarray*}
which implies \eqref{eq:var} by integrating over $s$.

Finally,
\begin{eqnarray}
\label{eq:covxdecay}
\quad \operatorname{cov}(X_h,X_0)&=&\operatorname
{cov}\biggl(\int_{M_d^-}\int
_{-\infty}^he^{A(h-s)}\Lambda(dA,ds),\int_{M_d^-}\int_{-\infty
}^0e^{-As}\Lambda(dA,ds) \biggr)\nonumber
\\
&=&\operatorname{cov}\biggl(\int_{M_d^-}\int_{-\infty
}^0e^{A(h-s)}\Lambda
(dA,ds),\int_{M_d^-}\int_{-\infty}^0e^{-As}\Lambda(dA,ds)
\biggr)\nonumber
\\[-8pt]\\[-8pt]
&=&\int_{M_d^-}e^{Ah}\biggl(\int_{-\infty}^0e^{-As}\biggl(\Sigma
+\int_{\mathbb{R}^d}xx^*\nu(dx)\biggr) e^{-A^*s}ds\biggr)\pi
(dA)\nonumber
\\
&=&-\int_{M_d^-}e^{Ah}(\mathscr{A}(A))^{-1}\biggl(\Sigma+\int
_{\mathbb{R}
^d}xx^*\nu(dx)\biggr)\pi(dA)\nonumber,
\end{eqnarray}
since $\Lambda$ is a L\'evy basis, and hence the random measures
$\Lambda$ on $M_d^-\times(0,h]$ and on $M_d^-\times(-\infty,0]$ are
independent.

From \eqref{eq:covxdecay} one obtains
\begin{eqnarray*}
&&\biggl\|\int_{M_d^-}e^{Ah}\biggl(\int_{-\infty}^0e^{-As}
\biggl(\Sigma+\int_{\mathbb{R}^d}xx^*\nu(dx)\biggr) e^{-A^*s}\,ds
\biggr)\pi
(dA)\biggr\|
\\
&&\qquad\leq\int_{M_d^-}\int_{-\infty}^0\kappa(A)^2e^{\rho
(A)(2s-h)}\,ds\,\pi(dA)\biggl\|\Sigma+\int_{\mathbb{R}^d}xx^*\nu
(dx)\biggr\|
\\
&&\qquad\leq\int_{M_d^-}\frac{\kappa(A)^2}{2\rho(A)}\pi
(dA)\biggl\|\Sigma+\int_{\mathbb{R}^d}xx^*\nu(dx)\biggr\|<\infty.
\end{eqnarray*}
Therefore $\lim_{h\to\infty}e^{Ah}=0$ for all $A\in M_d^-$ and
dominated convergence establish \eqref{eq:acovXzero}.
\end{pf}

\subsection{``SDE representation'' and some important path
properties}\label{sec:33}
In this section we show for a supOU process $X$ a representation which
generalises the SDE that governs OU-type processes, and we derive
important path properties of $X$. The ``SDE representation''---identity \eqref{eq:supousde}
below---has been conjectured in the
univariate case in \cite{barndorffnielsen01}, where neither a proof
nor conditions for its validity have been given. Below we are able to
show these results for finite variation L\'evy bases, which are
naturally appearing in applications like stochastic volatility
modeling. The properties which we establish are especially important
in the context of integration, since they imply that, if $X$ is the
integrator, then pathwise Lebesgue integration can be carried out, and,
when $X$ is the integrand, the theory of stochastic integrals of c\`
adl\`ag processes with respect to semimartingales (see \cite
{Protter2004}, for instance) respectively the $L^2$-theory of, for
example, \cite{Oksendal1998} applies. Likewise, the integrated process
is of importance in certain applications (see Section \ref{sec:finecon}).

Below the filtration $(\mathscr{F}_t)_{t\in\mathbb{R}}$ generated by
$\Lambda$ is defined by $\mathscr{F}_t$ being the $\sigma$-algebra
generated by the set of random variables $\{\Lambda(B)\dvtx B\in\mathscr
{B}(M_d^-\times(-\infty,t])\}$ for $t \in\mathbb{R}$.

\begin{Theorem}\label{th:pathprop}
Let $X$ be a supOU process as in Proposition \ref{th:supoufv}. Then:
\begin{longlist}[(iii)]
\item[(i)] $X_t(\omega)$ is $\mathscr{B}(\mathbb{R})\times\mathscr{F}$
measurable as a function of $t\in\mathbb{R}$ and $\omega\in\Omega$ and
adapted to the filtration $(\mathscr{F}_t)_{t\in\mathbb{R}}$
generated by
$\Lambda$.

\item[(ii)] If
%
\begin{equation}\label{eq:condbound}
\int_{M_d^-}\kappa(A)\pi(dA)<\infty,
\end{equation}
the paths of $X$ are locally uniformly bounded in $t$ for every $\omega
\in\Omega$.

Furthermore, $X_t^{+}=\int_0^tX_s\,ds$ exists for all $t\in\mathbb
{R}^+$ and
\begin{eqnarray}
\label{eq:intX} X_t^+&=&\int_{M_d^-}\int_{-\infty}^t
A^{-1}e^{A(t-s)}\Lambda(dA,ds)-\int_{M_d^-}\int_{-\infty}^0
A^{-1}e^{-As}\Lambda(dA,ds)\nonumber
\\[-8pt]\\[-8pt]
&&{}-\int_{M_d^-}\int_{0}^t A^{-1}\Lambda
(dA,ds)\nonumber.
\end{eqnarray}

\item[(iii)] Provided that
\begin{eqnarray} \label{eq:exZcond}
\int_{M_d^-}\frac{(\|A\| \vee1)\kappa(A)}{\rho(A)}\pi(dA)<\infty
\end{eqnarray}
and
\begin{eqnarray}\label{eq:boundZcond}
\int_{M_d^-}\|A\| \kappa(A)\pi(dA)<\infty
\end{eqnarray}
it holds that
%
\begin{equation}\label{eq:supousde}
X_t=X_0+\int_0^t Z_u \,du+L_t,
\end{equation}
where $L$ is the underlying L\'evy process and
%
\begin{equation}
Z_u=\int_{M_d^-}\int_{-\infty}^uAe^{A(u-s)}\Lambda(dA,ds)
\end{equation}
for all $u\in\mathbb{R}$ with the integral existing $\omega$-wise.

Moreover, the paths of $X$ are c\`adl\`ag and of finite variation on compacts.
\end{longlist}
\end{Theorem}

\begin{pf}
(i) is immediate from the definition of $X_t$ as a Lebesgue integral
and the measurability properties of the integrand $e^{A(t-s)}x
1_{\mathbb{R}
^+}(t-s)$ which as a function of $t,s,A,x$ is $\mathscr{B}(\mathbb
{R}\times
\mathbb{R}\times M_d^-\times\mathbb{R}^d)$-$\mathscr{B}(\mathbb
{R}^d)$-measurable.

(ii) We first show local uniform boundedness of $X$. Choose arbitrary
$T_1,T_2\in\mathbb{R}$ with $T_1<T_2$. Then
\begin{eqnarray*}f_{T_1,T_2}(A,s,x)&:=&\sup_{t\in[T_1,T_2]}\bigl\|
e^{A(t-s)}x1_{\mathbb{R}^+}(t-s)\bigr\|
\\
&\leq&\bigl(\kappa
(A)e^{-\rho
(A)(T_1-s)}1_{(-\infty,T_1]}(s)+\kappa(A)1_{(T_1,T_2]}(s)\bigr)\|x\|
\end{eqnarray*}
for all $A\in M_d^-(\mathbb{R}), s\in\mathbb{R}$ and $x\in\mathbb
{R}^d$ and
\begin{eqnarray*}
\sup_{t\in[T_1,T_2]}\|X_t\|&\leq&\int_{M_d^-}\int
_\mathbb{R}f_{T_1,T_2}(A,s,\gamma_0) \,ds\,\pi(dA)
\\
&&{}+\int
_{\mathbb{R}^d}\int
_{M_d^-}\int_\mathbb{R}f_{T_1,T_2}(A,s,x)\mu(dx,dA,ds).
\end{eqnarray*}
Therefore we only have to show the $\omega$-wise existence and
finiteness of the integral on the right-hand side. This is, however, an
immediate consequence of the above upper bound, \eqref{eq:condbound},
Proposition \ref{th:exlebsti} and arguments as in the proof of Theorem
\ref{th:supouex} noting that
\begin{eqnarray*}
&& \int_{M_d^-}\int_{T_1}^{T_2}\int_{\mathbb{R}^d}\bigl(1\wedge\kappa
(A)\|x\|
\bigr)\nu(dx)\,ds\,\pi(dA)
\\
&&\quad\quad\leq(T_2-T_1)\biggl(\int
_{M_d^-}\int_{\|x\|\leq1}\kappa(A)\|x\|\nu(dx)\pi(dA)
\\
&&{}\qquad\quad\hspace*{46pt}+\int
_{M_d^-}\int_{\|x\|> 1}1\nu(dx)\pi(dA)\biggr).
\end{eqnarray*}

Turning to $X_t^+$ the existence follows immediately from the local
boundedness. Noting that we have actually proved the local boundedness of
\[
\int_{M_d^-}\int_{-\infty}^t\bigl\|e^{A(t-s)}\gamma_0\bigr\| \,ds\,\pi(dA)+\int
_{\mathbb{R}^d}\int_{M_d^-}\int_{-\infty}^t\bigl\|e^{A(t-s)} x\bigr\|\mu(dx,dA,ds)
\]
above, we can use Fubini to obtain
\begin{eqnarray*}
X_t^+&=&\int_{M_d^-}\int_{-\infty}^t\int_{0\vee u}^te^{A(s-u)}\gamma
_0\,ds \,du\,\pi(dA)
\\
&&{}+\int_{\mathbb{R}^d}\int
_{M_d^-}\int
_{-\infty}^t\int_{0\vee u}^te^{A(s-u)} x\,ds\,\mu(dx,dA,du)
\\
&=&\int_{M_d^-}\int_{-\infty}^t A^{-1}e^{A(s-u)}\gamma_0
\big|_{s=(0\vee u)}^t \,du\,\pi(dA)
\\
&&{}+\int_{\mathbb{R}^d}\int_{M_d^-}\int
_{-\infty}^t A^{-1}e^{A(s-u)} x\big|_{s=(0\vee u)}^t\,ds\,\mu(dx,dA,du),
\end{eqnarray*}
which establishes \eqref{eq:intX} by straightforward
calculations.

(iii) Using similar calculations as before, the existence of $Z_u$ as
an $\omega$-wise integral follows from Proposition \ref{th:exlebsti}
and \eqref{eq:exZcond}. Similarly to (ii) one sees that under \eqref
{eq:boundZcond} $Z$ is locally uniformly bounded in $u$. Hence, one can
use Fubini to obtain
\begin{eqnarray*}
\int_0^t Z_u \,du&=&\int_{\mathbb{R}^d}\int_{M_d^-}\int_{-\infty
}^t\int
_{0\vee s}^tAe^{A(u-s)}x\,du\,\mu(dx,dA,ds)
\\
&&{}+\int_{M_d^-}\int_{-\infty
}^t\int_{0\vee s}^tAe^{A(u-s)}\gamma_0\,du\,ds\,\pi(dA)
\\
&=&\int_{\mathbb{R}^d}\int_{M_d^-}\int_{-\infty}^t
e^{A(u-s)}x\big|_{u=(0\vee s)}^t\mu(dx,dA,ds)
\\
&&{}+\int_{M_d^-}\int
_{-\infty}^t e^{A(u-s)}\gamma_0\big|_{u=(0\vee s)}^t\,ds\,\pi
(dA)=X_t-X_0-L_t
\end{eqnarray*}
which establishes \eqref{eq:supousde}. That $X$ has c\'adl\'ag paths
of finite variation is now an immediate consequence of this integral
representation.
\end{pf}

\begin{Remark}
(i) Condition \eqref{eq:condbound} is always true if $\pi$ is
concentrated on the normal matrices or on finitely many rays and hence
especially in dimension $d=1$. Moreover, it could be replaced by the
weaker but rather impracticable condition that
$f_{[T_1,T_2]}(A,s,x)\wedge1$ is integrable with respect to $\pi
\times\lambda\times\nu$ and that, for any fixed $x$,
$f_{[T_1,T_2]}(A,s,x)$ is integrable with respect to $\pi\times
\lambda$ (cf. \cite{marcusrosinski}, Proposition~2.1, for a very
related result whose proof is similar in spirit to ours, but uses a
series representation instead of the L\'evy--It\^o decomposition).

(ii) Intuitively \eqref{eq:boundZcond} means that $\pi$ does not
place too much mass on the elements of $M_d^-$ with high norm and thus
very fast exponential decay rates.

If $\pi$ is concentrated on the normal matrices or finitely many
diagonalizable rays, then \eqref{eq:exZcond} and \eqref
{eq:boundZcond} become
\begin{eqnarray}\label{eq:exZcondnormal}
\quad -\int_{M_d^{-}}\frac{(\|A\| \vee1)}{\max{\Re(\sigma(A))}}\pi
(dA)<\infty
\quad\mbox{and}\quad \int_{M_d^-}\|A\| \pi(dA)<\infty.
\end{eqnarray}
In particular, the second condition simply means that $\pi$ has a
finite first moment.

If $\pi$ is concentrated on $\mathbb{S}_d^{--}$, then we have $\|A\|
=-\min
(\sigma(A)),$ and \eqref{eq:exZcond} becomes
\begin{eqnarray} \label{eq:exZcondsym}
\int_{\mathbb{S}_d^{--}}\frac{(\min(\sigma(A)) \wedge-1)}{\max
{(\sigma
(A))}}\pi(dA)<\infty,
\end{eqnarray}
so it can be seen as a condition on the spread between the different
exponential decay rates measured by the eigenvalues. It is easy to see
that in dimension \mbox{$d=1$}, it is equivalent to $\int_{\mathbb
{R}^-}(-1/A)
\pi(dA)<\infty$, which is part of the necessary and sufficient
conditions for the existence of the supOU process.

(iii) It is very easy to construct examples when our sufficient
conditions for $X$ to exist as an $\omega$-wise integral are
satisfied, but neither \eqref{eq:exZcond} nor \eqref{eq:boundZcond}
for $Z$ above. Take, for example, $\pi$ concentrated on $v_n=
\left({ - n \atop 0}\enskip{ 0 \atop -1}
\right)
$ with $\pi(v_n)=6/(\pi^2n^2)$. In such a case we unfortunately do
not know whether $Z$ exists because our previously employed techniques
seem, at best, to give a necessary condition of the type \eqref
{eq:exZcond} involving $j(A)\wedge1$, and hence the necessary
conditions for $Z$ to exist would be implied by the sufficient ones for
$X$. Therefore, we also refrained from giving necessary conditions in
this section.
\end{Remark}

\subsection{Examples and long-range dependence}

Like in the univariate case, the expression \eqref{eq:covar} does not
imply that we necessarily have an exponential decay of the
autocovariance function and thus a short memory process. On the
contrary we can easily obtain a long memory process, as the following
examples exhibit. Note that this illustrates that \eqref{eq:acovXzero}
is not obvious and indeed requires a detailed proof as above.

Apart from showing that multivariate supOU processes may exhibit
long-range dependence, the purpose of this section is to analyze some
concrete examples and their properties.

Regarding long-range dependence, there is, unfortunately, basically no
general theory developed in the multivariate case so far. Below we mean
by long-range dependence (or long memory) simply that at least one
element of the autocovariance function decays asymptotically like
$h^{-\alpha}$ for the lag $h$ going to infinity and for some $\alpha
\in(0;1)$. Intuitively this should clearly be a case when one may
appropriately speak of long-range dependence. Establishing a general
theory for multivariate long-range dependence seems to be very
important, but is beyond the scope of this paper.

\begin{Example}\label{ex31}
 Let $\Lambda$ be a $d$-dimensional L\'evy basis with
generating quadruple $(\gamma,\Sigma,\nu,\pi)$ with $\nu$
satisfying $\int_{\mathbb{R}^d}\|x\|^2\nu(dx)<\infty$ and $\pi$ being
given as the distribution of $RB$ with a diagonalizable $B\in M_d^-$
and $ R$ being a real $\Gamma(\alpha,\beta)$-distributed random
variable with $\alpha>1,\beta\in\mathbb{R}^+\backslash\{0\}$.
Hence, $R$
has probability density $f(r)=\frac{\beta^\alpha}{\Gamma(\alpha
)}r^{\alpha-1}e^{-\beta r} 1_{\mathbb{R}^+}(r)$, and from
\begin{eqnarray*}
&& -\int_{M_d^{-}}\frac{1}{\max(\Re(\sigma(A)))}\pi(dA)
\\
&&\qquad =\frac
{-\beta^\alpha}{\max(\Re(\sigma(B)))\Gamma(\alpha)}\int
_{\mathbb{R}
^+}r^{\alpha-2}e^{-\beta r}\,dr
\\
&&\quad\quad=\frac{-\beta^\alpha
}{\max(\Re(\sigma(B)))\Gamma(\alpha)}\cdot\frac{\Gamma(\alpha
-1)}{\beta^{\alpha-1}}=\frac{-\beta}{\alpha\max(\Re(\sigma(B)))}
\end{eqnarray*}
we conclude that \eqref{eq:condmoms} holds. Hence, the process
$X_t=\int_{M_d^-}\int_{-\infty}^te^{A(t-s)}\Lambda(dA,\break ds)$ exists,
is stationary and has finite second moments. Similar calculations imply
that $\alpha>1$ is also necessary for $X_t$ to exist.

For the autocovariance function at positive lags $h$ we find
\begin{eqnarray*}
\operatorname{cov}(X_h,X_0)&=&-\int_{M_d^-}e^{Ah}(\mathscr
{A}(A))^{-1}\operatorname{vec}
\biggl(\Sigma+\int_{\mathbb{R}^d}xx^*\nu(dx)\biggr)\pi(dA)
\\
&=&\int_{\mathbb{R}^+}e^{Bhr-\beta I_d r}r^{\alpha-2}\,dr\,\biggl(-\frac
{\beta
^\alpha}{\Gamma(\alpha)}\mathscr{B}^{-1}\biggl(\Sigma+\int
_{\mathbb{R}
^d}xx^*\nu(dx)\biggr)\biggr)
\end{eqnarray*}
with $\mathscr{B}\dvtx M_d(\mathbb{R})\to M_d(\mathbb{R}), X\mapsto
BX+XB^*$. Let
now $U\in GL_d(\mathbb{C})$ and $\lambda_1,\lambda_2,\break
\ldots
,\lambda_d\in(-\infty,0)+i\mathbb{R}$ be such that
\[
UBU^{-1}=
\pmatrix{
\lambda_1 & 0 & \cdots&0\cr
0 &\lambda_2 & \ddots&\vdots\cr
\vdots& \ddots&\ddots& 0\cr
0 &\cdots& 0 &\lambda_d
}
.
\]
Then, from $\int_0^\infty t^{z-1}e^{-kt}\,dt=\Gamma(z)k^{-z}$ for all
$z,k\in(0,\infty)+i\mathbb{R}$, where the power is defined via the
principal branch of the complex logarithm (see \cite{AbramowitzStegun}, page~255), we obtain that
\begin{eqnarray*}
&& \int_{\mathbb{R}^+}e^{Bhr-\beta I_d r}r^{\alpha-2}\,dr
\\
&&\quad\quad
=U\int_{\mathbb{R}^+}\exp\left(-r\left(\beta I_d -
\pmatrix{
\lambda_1 & 0 & \cdots&0\cr
0 &\lambda_2 & \ddots&\vdots\cr
\vdots& \ddots&\ddots& 0\cr
0 &\cdots& 0 &\lambda_d
}
h\right)\right) r^{\alpha-2}\,dr\,U^{-1}
\\
&&\qquad=\Gamma(\alpha-1)U\left(\beta I_d -
\pmatrix{
\lambda_1 & 0 & \cdots&0\cr
0 &\lambda_2 & \ddots&\vdots\cr
\vdots& \ddots&\ddots& 0\cr
0 &\cdots& 0 &\lambda_d
}
h\right)^{1-\alpha}U^{-1}
\\
&&\quad\quad=\Gamma(\alpha-1)
(\beta I_d -Bh)^{1-\alpha}.
\end{eqnarray*}
Above the $(1-\alpha)$th power of a matrix is understood to be defined
via spectral calculus as usual.

Hence,
\[
\operatorname{cov}(X_h,X_0)=-\frac{\beta^\alpha}{\alpha-1}
(\beta I_d
-Bh)^{1-\alpha}\mathscr{B}^{-1}\biggl(\Sigma+\int_{\mathbb{R}
^d}xx^*\nu(dx)\biggr),
\]
and thus we have a polynomially decaying autocovariance function. For
$\alpha\in(1,2)$ we obviously get long memory.

Another question is whether we have the nice path properties of Theorem
\ref{th:pathprop}. Hence, assume additionally that $\int_{\|x\|\leq
1}\|x\|\nu(dx)<\infty$. In our example condition \eqref
{eq:condbound} is trivially satisfied and so the paths of $X$ are
locally uniformly bounded in~$t$. Regarding condition \eqref
{eq:exZcondnormal} the second part is equivalent to
\[
\|B\|\int_0^\infty r^\alpha e^{-\beta r}\,dr<\infty,
\]
which is always true, as any Gamma distribution has a finite mean.
Denoting the density of the $\Gamma(\alpha,\beta)$-distribution by
$f_{\alpha,\beta}(r)$, one obtains for the first part
\begin{eqnarray*}
-\int_0^\infty\frac{(r\|B\|\vee1)}{r\max(\Re(\sigma
(B)))}f_{\alpha,\beta}(r)\,dr&=&-\int_0^{\|B\|^{-1}} \frac{1}{r\max
(\Re(\sigma(B)))}f_{\alpha,\beta}(r)\,dr
\\
&&{}-\int_{\|B\|^{-1}}^\infty
\frac{\|B\|}{\max(\Re(\sigma(B)))}f_{\alpha,\beta}(r)\,dr,
\end{eqnarray*}
which is obviously finite. Hence, the conditions of Theorem \ref
{th:pathprop}(iii) are satisfied and thus the paths are c\`adl\`ag and
of finite variation, and \eqref{eq:supousde} is valid.
\end{Example}

\begin{Example}\label{ex32}
The previous example has an immediate extension to the case when $\pi$
is concentrated on several rays instead of a single one as above.
Assume we have $w_1,\ldots, w_m\in[0,1]$ with $\sum_{i=1}^m w_i=1$
and diagonalizable $B_1,\ldots, B_m\in M_d^-$, and define $\pi_i$ to
be the probability measure of the random variable $R_i B_i$ with $R_i$
being $\Gamma(\alpha_i,\beta_i)$ distributed with $\alpha_i>1$,
$\beta_i\in\mathbb{R}^+\backslash\{0\}$. If $\nu$ is as above and
$\pi
=\sum_{i=1}^m w_i\pi_i$, we get for the multivariate supOU process $X$
\[
\operatorname{cov}(X_h,X_0)=-\sum_{i=1}^m\biggl(\frac{w_i\beta
_i^{\alpha
_i}}{\alpha_i-1}(\beta_i I_d -B_ih)^{1-\alpha
_i}\mathscr{B}_i^{-1}\biggr)\biggl(\Sigma+\int_{\mathbb
{R}^d}xx^*\nu
(dx)\biggr)
\]
with $\mathscr{B}_i\dvtx M_d(\mathbb{R})\to M_d(\mathbb{R}), X\mapsto
B_iX+XB_i^*$.

Assuming now $\int_{\|x\|\leq1}\|x\|\nu(dx)<\infty$, it is likewise
straightforward to see that conditions \eqref{eq:condbound} and \eqref
{eq:exZcondnormal} are satisfied. Hence, the paths of $X$ are locally
uniformly bounded in $t$, c\`adl\`ag and of finite variation, and
\eqref{eq:supousde} is valid.
\end{Example}

\begin{Example}\label{ex33}
A similar result can be obtained if we restrict the mean reversion
parameter $A$ to the strictly negative definite matrices $\mathbb{S}_d^{--}$
and define $\pi$ as a probability distribution on the proper convex
cone $\mathbb{S}_d^{--}$ as follows. Let $\mathbf{S}_d^{--}$ denote the
intersection of the unit sphere in $\mathbb{S}_d$ with $\mathbb
{S}_d^{--}$, let
$\alpha\dvtx \mathbf{S}_d^{--}\to(1,\infty)$, $\beta\dvtx \mathbf
{S}_d^{--}\to(0,\infty)$ be measurable mappings and $w$ a probability
distribution on $\mathbf{S}_d^{--}$ such that
%
\begin{equation}\label{eq:ex3excond}
-\int_{\mathbf{S}_d^{--}}\frac{\beta(v)}{\alpha(v)\max(\sigma
(v))}w(dv)<\infty.
\end{equation}
Now define $\pi$ via
\[
\pi(B)=\int_{\mathbf{S}_d^{--}}\int_0^\infty1_{B}(rv)\frac{\beta
(v)^{\alpha(v)}}{\Gamma(\alpha(v))}r^{\alpha(v)-1}e^{-\beta(v)r}\,dr\,w(dv)
\]
for any Borel set $B\in M_d^-(\mathbb{R})$. Then $\pi$ is a probability
distribution concentrated on~$\mathbb{S}_d^{--}$.

Using this $\pi$ in the above set-up means that the mean reversion
parameter is no longer necessarily restricted to finitely many rays.
Moreover, similar calculations to the ones in Example \ref{ex31} give
\begin{eqnarray*}
-\int_{M_d^{-}}\frac{1}{\max(\Re(\sigma(A)))}\pi(dA)=\int
_{\mathbf{S}_d^{--}}\frac{-\beta(v)}{\alpha(v)\max(\sigma
(v))}w(dv)<\infty.
\end{eqnarray*}
Hence, \eqref{eq:condmoms} holds and the process $X_t=\int
_{M_d^-}\int_{-\infty}^te^{A(t-s)}\Lambda(dA,ds)$ exists, is
stationary and has finite second moments. Likewise we get for the
autocovariance function
\begin{eqnarray*}
\operatorname{cov}(X_h,X_0)&=&-\biggl(\int_{\mathbf{S}_d^{--}}\frac
{\beta
(v)^{\alpha(v)}}{\alpha(v)-1}\biggl(\beta(v) I_d -vh
\biggr)^{1-\alpha(v)}(\mathscr{V}(v))^{-1}w(dv)\biggr)
\\
&&{} \hspace*{6pt}\times
\biggl(\Sigma+\int_{\mathbb{R}^d}xx^*\nu(dx)\biggr)
\end{eqnarray*}
with $\mathscr{V}(v)\dvtx M_d(\mathbb{R})\to M_d(\mathbb{R}), X\mapsto vX+Xv^*$.

Turning to the path properties, assume now that $\int_{\|x\|\leq1}\|
x\|\nu(dx)<\infty$. Again condition \eqref{eq:condbound} is
trivially satisfied and so the paths of $X$ are locally uniformly
bounded in $t$. Regarding condition \eqref{eq:exZcondnormal} the
second part becomes
\[
\int_{M_d^-}\|A\| \pi(dA)=\int_{\mathbf{S}_d^{--}}\int_0^\infty
rf_{\alpha(v),\beta(v)}(r)\,dr\,w(dv)=\int_{\mathbf{S}_d^{--}}\frac
{\alpha(v)}{\beta(v)}w(dv),
\]
and for the first part one obtains
\begin{eqnarray}\label{eq:exzexample}
&&-\int_{M_d^{-}}\frac{(\|A\| \vee1)}{\max{\Re(\sigma(A))}}\pi\nonumber
(dA)
\\
&&\qquad =-\int_{\mathbf{S}_d^{--}}\int_0^{1} \frac{1}{r\max(\sigma
(v))}f_{\alpha(v),\beta(v)}(r)\,dr\,w(dv)
\\
&&\qquad\quad {}-\int_{\mathbf
{S}_d^{--}}\int_{1}^\infty\frac{1}{\max(\sigma(v))}f_{\alpha
(v),\beta(v)}(r)\,dr\,w(dv).\nonumber
\end{eqnarray}
The first summand is finite due to \eqref{eq:ex3excond} and the second
one is finite if the integral $-\int_{\mathbf{S}_d^{--}}(1/\max
(\sigma(v)))w(dv)$ is finite.
Hence, provided
\begin{eqnarray*}
-\int_{\mathbf{S}_d^{--}}\frac{1}{\max(\sigma(v))}w(dv)<\infty
\quad\mbox{and}\quad\int_{\mathbf{S}_d^{--}}\frac
{\alpha(v)}{\beta(v)}w(dv)<\infty,
\end{eqnarray*}
the conditions of Theorem \ref{th:pathprop}(iii) are satisfied, and
thus the paths are c\`adl\`ag and of finite variation and \eqref
{eq:supousde} is valid.

Based on this we can easily give an example where we know that the
supOU process exists due to Proposition \ref{th:supoufv}, but the
conditions of Theorem \ref{th:pathprop}(iii) are not satisfied.
Assume $w$ is a discrete distribution concentrated on the points
\[
v_n=
\pmatrix{ -1 & 0 &0\cr
 0 & -1+(3n)^{-1} &0\cr
  0 & 0 & -1/2
}
,\qquad n\in\mathbb{N},
\]
and that $w(v_n)=\frac{6}{\pi^2 n^2 }$, $\alpha(v_n)=2$ and $\beta
(v_n)=n^{-1}$. Then we have that
\[
-\int_{\mathbf{S}_d^{--}}\frac{\beta(v)}{\alpha(v)\max(\sigma
(v))}w(dv)=\frac{6}{\pi^2}\sum_{n=1}^\infty n^{-3}<\infty,
\]
but
\[
\int_{M_d^-}\|A\| \pi(dA)=\frac{12}{\pi^2}\sum_{n=1}^\infty
n^{-1}=\infty,
\]
and hence condition \eqref{eq:boundZcond} is not satisfied. Observe
that this means that the probability measure $\pi$ we have constructed
does not have a first moment, although it is defined via a polar
representation where the radial parts are all univariate Gamma
distributions. Moreover, it is easy to see that \eqref{eq:exzexample}
is finite and so $Z$ exists and thus it is only the local uniform
boundedness our sufficient conditions fail to provide when trying to
show \eqref{eq:supousde}. Showing that \eqref{eq:supousde} is indeed
not valid seems to be a very delicate issue, as already remarked.
\end{Example}

\begin{Example}\label{ex34}
Let $\Lambda$ be now a two-dimensional L\'evy basis with generating
quadruple $(\gamma,\Sigma,\nu,\pi)$ with $\nu$ satisfying $\int
_{\mathbb{R}^2}\|x\|^2\nu(dx)<\infty$. We restrict the mean reversion
parameter $A$ to $\mathbb{D}^{--}_2$, the $2 \times 2$ diagonal matrices with
strictly negative entries on the diagonal. Hence, $\pi$ is a measure
on $\mathbb{D}^{--}_2$, which can be identified with $(\mathbb{R}^{--})^2$,
and we assume that $\pi$ has Lebesgue density\vspace*{-1pt}
\begin{eqnarray*}
&&\pi(da_1,da_2)
\\[-2pt]
&&\qquad =\frac{\beta_1^{\alpha_1}\beta_2^{\alpha_2}}{\Gamma
(\alpha_1)\Gamma(\alpha_2)}(-a_1)^{\alpha_1-1}(-a_2)^{\alpha
_2-1}e^{\beta_1a_1+\beta_2a_2}1_{(\mathbb{R}^{--})^2}(a_1,a_2)\,da_1 \,da_2
\end{eqnarray*}
with $\alpha_1,\alpha_2>1$ and $\beta_1,\beta_2>0$. So the diagonal
elements are independent, and their absolute values follow Gamma distributions.
We obtain
\begin{eqnarray*}
&& -\int_{\mathbb{D}_2^{--}}\frac{1}{\max(\Re(\sigma(A)))}\pi(dA)
\\[-2pt]
&&\quad\quad=\int_0^\infty\int_0^\infty\frac{1}{\min
(a_1,a_2)}\frac{\beta_1^{\alpha_1}\beta_2^{\alpha_2}}{\Gamma
(\alpha_1)\Gamma(\alpha_2)}(a_1)^{\alpha_1-1}(a_2)^{\alpha
_2-1}e^{-\beta_1a_1-\beta_2a_2}\,da_1\,da_2
\\[-2pt]
&&\quad\quad\leq\int_0^\infty\frac{\beta_1^{\alpha_1}}{\Gamma
(\alpha_1)}(a_1)^{\alpha_1-2}e^{-\beta_1a_1}\,da_1\int_0^\infty\frac
{\beta_2^{\alpha_2}}{\Gamma(\alpha_2)}(a_2)^{\alpha_2-1}e^{-\beta
_2a_2}\,da_2
\\[-2pt]
&&{}\quad\quad\quad+\int_0^\infty\frac{\beta_1^{\alpha_1}}{\Gamma
(\alpha_1)}(a_1)^{\alpha_1-1}e^{-\beta_1a_1}\,da_1\int_0^\infty\frac
{\beta_2^{\alpha_2}}{\Gamma(\alpha_2)}(a_2)^{\alpha_2-2}e^{-\beta
_2a_2}\,da_2<\infty.
\end{eqnarray*}
Hence, \eqref{eq:condmoms} holds, and the process $X_t=\int
_{M_d^-}\int_{-\infty}^te^{A(t-s)}\Lambda(dA,ds)$ exists, is
stationary and has finite second moments.

Let us now consider the individual components $X_{1,t}$, $X_{2,t}$ of
$X_t$. Denote by $P_1\dvtx \mathbb{R}^2\to\mathbb{R}, (x_1,x_2)^*\mapsto
x_1$ the
projection onto the first coordinate and define an $\mathbb{R}$-valued
L\'
evy basis $\Lambda_1$ on $\mathbb{R}^{--}\times\mathbb{R}$ via
$\Lambda
_1(da_1,ds)=P(\Lambda(P_1^{-1}(da_1),ds)$ and a L\'evy measure $\nu
_1$ on $\mathbb{R}$ via $\nu_1(dx_1)=\nu(P_1^{-1}(dx_1))$. Then
$\Lambda
_1$ has characteristic quadruple $(\gamma_1,\Sigma_{11},\nu_1,\pi
_1)$ with $\pi_1$ having Lebesgue density
\[
\pi_1(da_1)=\frac{\beta_1^{\alpha_1}}{\Gamma(\alpha
_1))}(-a_1)^{\alpha_1-1}e^{\beta_1a_1}1_{(\mathbb{R}^{--})}(a_1)\,da_1
\]
and
\[
X_{1,t}=\int_{\mathbb{R}^{--}}\int_{-\infty}^te^{a_1(t-s)}\Lambda
_1(da_1,ds).
\]
For the autocovariance function of the first component we get
\[
\operatorname{cov}(X_{1,h},X_{1,0})=\frac{\beta_1^{\alpha
_1}}{2(\alpha
_1-1)}(\beta_1 +h)^{1-\alpha_1}\biggl(\Sigma_{11}+\int
_{\mathbb{R}}x_1^2\nu_1(dx_1)\biggr),\qquad   h\in\mathbb{R}^+.
\]
An analogous result holds for the second component $X_{2,t}$ and we
have long memory in both components provided $\alpha_1,\alpha_2\in(1,2)$.

The importance of this example is, however, that we can model the
stationary distributions of $X_1$ and $X_2$, that is, the margins of
the stationary distribution of $X_t$, very explicitly by specifying the
margins of $\nu$, that is, $\nu_1$ and $\nu_2$. From \cite{barndorffnielsen01}, Theorem~3.1, Corollary~3.1, and
\cite{FasenetCklu2007}, Remark~2.2, we know that it is exactly all nondegenerate
self-decomposable distributions on $\mathbb{R}$ which arise as the
stationary distributions of the components. Moreover, these authors
provide formulae to calculate $\nu_1$ (or $\nu_2$) if one wants to
obtain a given stationary distribution for the component (alternatively
\cite{BarndorffetJensenetSoerensen1998}, Lemma 5.1, or the refinement
\cite{PigorschetStelzer2007b}, Theorem 4.9, can be used). Hence, one
can specify a two-dimensional supOU process with prescribed stationary
distributions of the components by calculating the required $\nu_1$
and $\nu_2$ and choosing $\nu$ accordingly. The easiest way to get a
possible $\nu$ is by specifying $\nu(dx_1,dx_2)=\nu_1(dx_1)\times
\delta_0(x_2)+\delta_0(x_1)\times\nu_2(dx_2)$ with $\delta_0$
denoting the Dirac distribution with unit mass at zero. In this case
the components of $X$ are independent. An easy way to get an
appropriate $\nu$ and allowing for dependence is to combine $\nu_1$
and $\nu_2$ using a L\'evy copula (see \cite{BarndorffLindner2007,KallsenTankov2006}).

Likewise it is again interesting to look at the path properties of
Theorem \ref{th:pathprop}. Assuming again $\int_{\|x\|\leq1}\|x\|\nu
(dx)<\infty$, condition \eqref{eq:condbound} is trivially satisfied,
and so the paths of $X$ are locally uniformly bounded in $t$. Regarding
the second part of condition \eqref{eq:exZcondnormal} we have that
\[
\int_{M_d^{-}}\|A\|\pi(dA)<\infty
\]
is equivalent to
\begin{eqnarray*}
&&\int_{(\mathbb{R}^+)^2}\max(a_1,a_2)f_{\alpha_1,\beta
_1}(a_1)f_{\alpha
_2,\beta_2}(a_2)\,da_1\,da_2
\\
&&\qquad \leq \int_{\mathbb{R}^+}a_1f_{\alpha
_1,\beta
_1}(a_1)\,da_1
+\int_{\mathbb{R}^+}a_2f_{\alpha_2,\beta
_2}(a_2)\,da_2<\infty,
\end{eqnarray*}
which is always true. Turning to \eqref{eq:exZcondsym}, it is implied by
\begin{eqnarray*}
&&\int_{(\mathbb{R}^+)^2}\frac{\max(a_1,a_2)}{\min
(a_1,a_2)}f_{\alpha
_1,\beta_1}(a_1)f_{\alpha_2,\beta_2}(a_2)\,da_1\,da_2
\\
&&\quad\quad
\leq\int_{(\mathbb{R}^+)^2}\frac{a_1+a_2}{a_1}f_{\alpha_1,\beta
_1}(a_1)f_{\alpha_2,\beta_2}(a_2)\,da_1\,da_2
\\
&&{}\quad\quad\quad+\int
_{(\mathbb{R}^+)^2}\frac{a_1+a_2}{a_2}f_{\alpha_1,\beta
_1}(a_1)f_{\alpha
_2,\beta_2}(a_2)\,da_1\,da_2<\infty,
\end{eqnarray*}
which is easily seen to be always true. Hence, the conditions of
Theorem \ref{th:pathprop}(iii) are satisfied, and thus the paths are
c\`adl\`ag and of finite variation and \eqref{eq:supousde} is valid.

Obviously this example has a straightforward extension to general
dimension $d$.
\end{Example}

\begin{Example}\label{ex5}
So far we have only studied cases where we could use Proposition \ref
{th:exprop} and did especially never have to bother with $\kappa(A)$
in the conditions of Theorem \ref{th:supouex}.

In this example we will present a case where the behavior of $\kappa
(A)$ is crucial and where we show how $\kappa$ and $\rho$ can be
specified in a measurable way.
We define the following sets:
\begin{eqnarray*}
\mathscr{D}^-_d&=&\{X\in M_d(\mathbb{R})\dvtx X \mbox{ is diagonal; all
diagonal elements are strictly negative,}
\\
&&\hspace*{18pt}\mbox{ pairwise distinct
and ordered such that } x_{ii}< x_{jj}\  \forall1\leq i\leq j\leq d \};
\\
\mathscr{S}_d&=&\{X\in GL_d(\mathbb{R})\dvtx \mbox{ the first nonzero element
in each column is } 1\};
\\
\mathscr{M}^-_d&=&\{SDS^{-1}\dvtx   S\in\mathscr{S}_d, D\in\mathscr
{D}^-_d\}.
\end{eqnarray*}
If $A=SDS^{-1}$ is in $\mathscr{M}_d^-$, the matrix $D$ consists of
the eigenvalues of $A$, and the columns of $S$ are the eigenvectors of
$A$. In principle there are many possible $S$ and $D$ if we only demand
$A=SDS^{-1}$. However, if we restrict ourselves to $S\in\mathscr
{S}_d, D\in\mathscr{D}^-_d$, then $S,D$ are unique, as elementary
linear algebra shows. This means that the map
\[
\mathfrak{M}\dvtx \mathscr{S}_d\times\mathscr{D}^-_d\to\mathscr
{M}_d^-, (S,D)\mapsto SDS^{-1}
\]
is bijective (and obviously continuous). We denote by $\mathfrak
{M}^{-1}=(\mathfrak{S},\mathfrak{D})$ the inverse mapping. Since
computing eigenvectors and eigenvalues are measurable procedures as are
the orderings and normalizations involved in obtaining the diagonal
matrix in $\mathscr{D}^-_d$ and the eigenvector matrix in $\mathscr
{S}_d$, all these mappings are measurable. Note also that $\mathscr
{D}^-_d$, $\mathscr{S}_d$, $\mathscr{M}^-_d$ are Borel sets.

Defining $\kappa\dvtx \mathscr{M}^-_d\to[1,\infty), A\mapsto\|
\mathfrak{S}(A)\|\|(\mathfrak{S}(A))^{-1}\|$, $\rho(A)=\break -\max(\Re
(\sigma(A)))$ gives, therefore, measurable mappings on $\mathscr
{M}^-_d$ satisfying $\|e^{As}\|\leq\kappa(A)e^{-\rho(A)s}$. Using
these definitions for $\kappa$ and $\rho$ one could now specify
probability distributions $\pi$ on $\mathscr{M}^-_d$ and check
whether condition \eqref{c2} is satisfied, and the associated supOU
process therefore exists.

However, in concrete situations it seems easier to specify a Borel
probability measure $\pi_{\mathscr{S}_d\times\mathscr{D}^-_d}$ on
$\mathscr{S}_d\times\mathscr{D}^-_d$ and define $\pi$ as its image
under $\mathfrak{M}$, i.e. $\pi(B)=\pi_{\mathscr{S}_d\times
\mathscr{D}^-_d}(\mathfrak{M}^{-1}(B))$ for all Borel sets $B$.
Assume $\pi_{\mathscr{S}_d\times\mathscr{D}^-_d}=\pi_{\mathscr
{S}_d}\times\pi_{\mathscr{D}^-_d}$ is the product of two probability
measures $\pi_{\mathscr{S}_d}$ on $\mathscr{S}_d$ and $\pi
_{\mathscr{D}^-_d}$ on $\mathscr{D}^-_d$.
Then we have
\begin{eqnarray*}
\int_{M_d^-}\frac{\kappa(A)^2}{\rho(A)}\pi(dA)<\infty
&\quad\Longleftrightarrow\quad& \int_{\mathscr{S}_d}\|S\|^2\|S^{-1}\|^2\pi
_{\mathscr{S}_d}(dS)<\infty\quad\mbox{and }
\\
&&{}-\int_{\mathscr
{D}^-_d}\frac{1}{\max(\Re(\sigma(D)))}\pi_{\mathscr
{D}^-_d}(dD)<\infty.
\end{eqnarray*}

That $\int_{\mathscr{S}_d}\|S\|^2\|S^{-1}\|^2\pi_{\mathscr
{S}_d}(dS)$ can be finite or infinite depending on the choice of $\pi
_{\mathscr{S}_d}$ is exhibited by the following example. Let $\pi
_{\mathscr{S}_2}$ be a discrete measure concentrated on the points
\[
S_n=
\pmatrix{ 1 &0 \cr n & 1
}
\quad\mbox{and}\quad p_n:=\pi_{\mathscr{S}_2}(S_n)=C_\alpha n^{-\alpha}
\qquad \forall n\in\mathbb{N}
\]
with $\alpha>1$ and $C_\alpha=1/\sum_{n=1}^\infty n^{-\alpha}$. Then
\[
S^{-1}_n=
\pmatrix{ 1 &0 \cr -n & 1
}
.
\]
Using the equivalence of all norms we get that
\[
\int_{\mathscr{S}_2}\|S\|^2\|S^{-1}\|^2\pi_{\mathscr
{S}_2}(dS)<\infty \quad\Longleftrightarrow\quad  C_\alpha\sum_{n=1}^\infty n^4
p_n<\infty \quad\Longleftrightarrow \quad \alpha>5.
\]

Returning to the general example with $\pi$ given via $\pi_{\mathscr
{S}_d}\times\pi_{\mathscr{D}^-_d}$ and turning to path properties,
we assume again $\int_{\|x\|\leq1}\|x\|\nu(dx)<\infty$. In this
finite variation case the existence conditions \eqref{fvc2} become
\[
\int_{\mathscr{S}_d}\|S\|\|S^{-1}\|\pi_{\mathscr{S}_d}(dS)<\infty
\quad\mbox{and}\quad -\int_{\mathscr{D}^-_d}\frac{1}{\max(\Re(\sigma
(D)))}\pi_{\mathscr{D}^-_d}(dD)<\infty.
\]
Furthermore, condition \eqref{eq:condbound} is always satisfied when
the existence conditions are satisfied and so the paths of $X$ are
locally uniformly bounded in $t$. Straightforward arguments show that
the conditions of Theorem \ref{th:pathprop}(iii) are satisfied, and
thus the paths are c\`adl\`ag and of finite variation and \eqref
{eq:supousde} is valid if
\begin{eqnarray*}
\int_{\mathscr{S}_d}\|S\|^2\|S^{-1}\|^2\pi_{\mathscr
{S}_d}(dS)<\infty,\qquad  -\int_{\mathscr{D}^-_d}\frac{\|D\|}{\max(\Re
(\sigma(D)))}\pi_{\mathscr{D}^-_d}(dD)<\infty
\end{eqnarray*}
 and
 \begin{eqnarray*}
 \int
_{\mathscr{D}^-_d}\|D\|\pi_{\mathscr{D}^-_d}(dD)<\infty.
\end{eqnarray*}
By polarly decomposing $\pi_{\mathscr{D}^-_d}$ into a measure on the
unit sphere in the diagonal matrices and a radial part, the long memory
specifications of the foregoing examples have straightforward
extensions to this set-up.
\end{Example}

\section{Positive semi-definite supOU processes}\label{sec:4}

Based on the previous section we now consider supOU processes which are
positive semi-definite at all times. The importance of such processes
is that they can be used to describe the random evolution of a latent
covariance matrix over time and, hence, they can be used in
multivariate models for heteroskedastic data, for example, the
stochastic volatility model of \cite{BarndorffetStelzer2009sv}.

Let us briefly recall that a $d\times d$ positive semi-definite OU-type
process (see \cite{BarndorffetStelzer2006}) is defined as the unique
c\`adl\`ag solution of the SDE
\[
d\Sigma_t=(A\Sigma_t+\Sigma_tA^*)\,dt+dL_t, \qquad \Sigma_0\in\mathbb{S}_d^+,
\]
with $A\in M_d(\mathbb{R})$ and $L$ being a $d\times d$ matrix subordinator
(see \cite{BarndorffetPerez2005}), that is, a L\'evy process in
$\mathbb{S}
_d$ with $L_t-L_s\in\mathbb{S}_d^+$ $\forall  s,t\in\mathbb
{R}^+,  s<t$. If
$E(\ln(\max(\|L_1\|,1)))<\infty$ and $\max(\Re(\sigma(A)))<0$,
the above SDE has the unique stationary solution
\[
\Sigma_t=\int_{-\infty}^te^{A(t-s)}\,dL_se^{A^*(t-s)}.
\]
That the linear operators $\mathbb{S}_d\to\mathbb{S}_d$ of the form
$Z\mapsto
AZ+ZA^*$ with some $A\in M_d(\mathbb{R})$ are {the} ones to be used
for positive semi-definite OU-type processes has been established in
\cite{PigorschetStelzer2007b}.

As just recalled, one has to restrict the driving L\'evy process to
matrix subordinators in order to obtain OU-type processes taking values
in the positive semi-definite matrices. Below we need to impose a
comparable condition on the L\'{e}vy basis to get positive
semi-definite supOU processes. Note that for a $d\times d$
matrix-valued L\'evy basis $\Lambda$ we denote by $\operatorname
{vec}(\Lambda)$
the $\mathbb{R}^{d^2}$-valued L\'evy basis given by $\operatorname
{vec}(\Lambda)(B)=\operatorname{vec}
(\Lambda(B))$ for all Borel sets $B$. Moreover, observe that
$\operatorname{tr}
(XY^*)$ [with $X,Y\in
M_d(\mathbb{R})$ and $\operatorname{tr}$ denoting the usual trace
functional] defines a
scalar product on $M_d(\mathbb{R})$ and that the $\operatorname{vec}$
operator is a
Hilbert space isometry between $M_d(\mathbb{R})$ equipped with this scalar
product and $\mathbb{R}^{d^2}$ with the usual Euclidean scalar product.

Positive semi-definite supOU processes are defined as processes of the
form \eqref{eq:psdsupou} below, which is the analogue of \eqref{eq:defsupou}.

\begin{Theorem}\label{th:exsupoupsd}
Let $\Lambda$ be an $\mathbb{S}_d$-valued L\'evy basis on
$M_d^-\times\mathbb{R}
$ with generating quadruple $(\gamma, 0,\nu,\pi)$ with $
\gamma_0:=\gamma-\int_{\|x\|\leq1} x\nu(dx)\
\in\mathbb{S}_d^+$ and $\nu$ being a L\'evy measure on $\mathbb{S}_d$
satisfying $\nu(\mathbb{S}_d\backslash\mathbb{S}_d^+)=0$,
\begin{equation}
\int_{\|x\|>1}\ln(\|x\|)\nu(dx)<\infty \quad\mbox{and}\quad   \int_{\|
x\|\leq1}\|x\|\nu(dx)<\infty.
\end{equation}
Moreover, assume there exist measurable functions $\rho\dvtx M_d^-\to
\mathbb{R}
^+\backslash\{0\}$ and $\kappa\dvtx M_d^-\to[1,\infty)$ such that
\begin{eqnarray}
&\|e^{As}\|\leq\kappa(A)e^{-\rho(A)s} \qquad \forall s\in\mathbb
{R}^+,
\pi\mbox{-almost surely},
\end{eqnarray}
 and
\begin{eqnarray}\label{eq:excondpsd}
 \int_{M_d^-}\frac{\kappa(A)^2}{\rho(A)}\pi
(dA)<\infty.
\end{eqnarray}

Then the process $(\Sigma_t)_{t\in\mathbb{R}}$ given by
\begin{eqnarray}\label{eq:psdsupou}
\Sigma_t&=&\int_{M_d^-}\int_{-\infty}^te^{A(t-s)}\Lambda
(dA,ds)e^{A^*(t-s)}\nonumber
\\
&=&\int_{M_d^-}\int_{-\infty}^te^{A(t-s)}\gamma
_0e^{A^*(t-s)} \,ds\,\pi(dA)
\\
&&{}+\int_{\mathbb{S}_d}\int_{M_d^-}\int_{-\infty}^te^{A(t-s)}
xe^{A^*(t-s)}\mu(dx,dA,ds)\nonumber
\end{eqnarray}
is well defined as a Lebesgue integral for all $t\in\mathbb{R}$ and
$\omega
\in\Omega$ and $\Sigma$ is stationary.

Moreover,
%
\begin{equation}\label{eq:psdsupouvec}
\operatorname{vec}(\Sigma_t)=\int_{M_d^-}\int_{-\infty
}^te^{(A\otimes
I_d+I_d\otimes A)(t-s)}\operatorname{vec}(\Lambda)(dA,ds),
\end{equation}
$\Sigma_t\in\mathbb{S}_d^+$ for all $t\in\mathbb{R}$ and the
distribution of
$\Sigma_t$ is infinitely divisible with characteristic function
\[
E(\exp(i\operatorname{tr}(u\Sigma_t))
)=\exp\biggl(i\operatorname{tr}
(u\gamma_{\Sigma,0})+\int_{\mathbb{S}_d}\bigl(e^{i\operatorname
{tr}(ux)}-1\bigr)\nu
_{\Sigma}(dx)\biggr), \qquad u\in\mathbb{S}_d,
\]
where
\begin{eqnarray}
\gamma_{\Sigma,0}&=&\int_{M_d^-}\int_0^\infty e^{As}\gamma
_0e^{A^*s}\,ds\,\pi(dA),
\\
\nu_{\Sigma}(B)&=&\int_{M_d^-}\int_{0}^\infty\int_{\mathbb{S}
_d^+}1_B(e^{As}xe^{A^*s})\nu(dx)\,ds\,\pi(dA)
\\
&&\eqntext{\mbox{for all Borel sets
} B\subseteq\mathbb{S}_d.}
\end{eqnarray}
\end{Theorem}

\begin{pf}
The equivalence of \eqref{eq:psdsupouvec} and \eqref{eq:psdsupou}
follows from standard results on the vectorization operator and the
tensor product (see \cite{Hornetal1991}).

Next we note that $e^{(A\otimes I_d+I_d\otimes A)(t-s)}=e^A\otimes e^A$
and that $\|e^{A}\otimes e^A\|=\|e^A\|^2$ (using the operator norm
associated with the Euclidean norm). Hence, all assertions except
$\Sigma_t\in\mathbb{S}_d^+$ for all $t\in\mathbb{R}^+$ follow
immediately from
Propositions \ref{th:exlebsti} and~\ref{th:supoufv}.

However, $\Sigma_t\in\mathbb{S}_d^+$ for all $t\in\mathbb{R}^+$ is now
immediate, since the integral exists $\omega$-wise,
$e^{As}Xe^{A^*s}\in\mathbb{S}_d^+$ $\forall A\in M_d(\mathbb
{R}),X\in\mathbb{S}
_d^+, s\in\mathbb{R}$ and $\mathbb{S}_d^+$ is a closed convex\break cone.
\end{pf}

\begin{Remark}
(i) As in Proposition \ref{th:exprop}, $\kappa(A)$ can be replaced by
$1$ and $\rho(A)$ by $-\max(\Re(\sigma(A)))$ in \eqref
{eq:excondpsd} [and also in \eqref{eq:excondmoms2psd} and \eqref
{eq:excondmompsd} below] provided $\pi$ is concentrated on the normal
matrices or finitely many diagonalizable rays.

(ii) Throughout this section we refrain from stating necessary
conditions, as they can be immediately inferred from the foregoing
sections and the arguments presented for the sufficient conditions.
\end{Remark}

Most importantly in the context of stochastic volatility models, which
involve stochastic integrals with $\Sigma$ as integrand, Theorem \ref
{th:pathprop} also has an analogue for positive semi-definite supOU processes.

\begin{Theorem}\label{th:pathproppd}
Let $\Sigma$ be the positive semi-definite supOU process of Theorem~\ref{th:exsupoupsd}. Then:\vspace*{-6pt}
\begin{longlist}[(iii)]
\item[(i)] $\Sigma_t(\omega)$ is $\mathscr{B}(\mathbb{R})\times\mathscr{F}$
measurable as a function of $t\in\mathbb{R}$ and $\omega\in\Omega$ and
adapted to the filtration $(\mathscr{F}_t)_{t\in\mathbb{R}}$
generated by
$\Lambda$.

\item[(ii)] If
%
\begin{equation}\label{eq:condboundpd}
\int_{M_d^-}\kappa(A)^2\pi(dA)<\infty,
\end{equation}
the paths of $\Sigma$ are locally uniformly bounded in $t$ for every
$\omega\in\Omega$.

Furthermore, $\Sigma_t^{+}=\int_0^t\Sigma_s\, ds$ exists for all $t\in
\mathbb{R}^+$ and
\begin{eqnarray}\label{eq:intSigma}
\Sigma_t^+&=&\int_{M_d^-}\int_{-\infty}^t (\mathbf{A}(A))^{-1}
\bigl(e^{A(t-s)}\Lambda(dA,ds)e^{A^*(t-s)}\bigr)\nonumber
\\
&&{}-\int
_{M_d^-}\int_{-\infty}^0 (\mathbf{A}(A))^{-1}(e^{-As}\Lambda
(dA,ds)e^{-A^*s})
\\
&&{}-\int_{M_d^-}\int_{0}^t (\mathbf
{A}(A))^{-1}\Lambda(dA,ds)\nonumber
\end{eqnarray}
with $\mathbf{A}(A)\dvtx \mathbb{S}_d\to\mathbb{S}_d, X\mapsto AX+XA^*$.

\item[(iii)] Provided that
\begin{eqnarray} \label{eq:exZcondpd}
-\int_{M_d^-}\frac{(\|A\| \vee1)\kappa(A)^2}{\rho(A)}\pi
(dA)&<\infty
\end{eqnarray}
and
\begin{eqnarray}\label{eq:boundZcondpd}
 \int_{M_d^-}\|A\| \kappa(A)^2\pi(dA)&<\infty
\end{eqnarray}
it holds that
%
\begin{equation}\label{eq:supousdepd}
\Sigma_t=\Sigma_0+\int_0^t Z_u du+L_t,
\end{equation}
where $L$ is the underlying matrix subordinator and
%
\begin{eqnarray}
Z_u&=&\int_{M_d^-}\int_{-\infty}^u\bigl(Ae^{A(u-s)}\Lambda
(dA,ds)e^{A^*(u-s)}\nonumber
\\[-8pt]\\[-8pt]
&&{}\hspace*{43pt}+e^{A(u-s)}\Lambda(dA,ds)e^{A^*(u-s)}A^*\bigr)\nonumber
\end{eqnarray}
for all $u\in\mathbb{R}$ with the integral existing $\omega$-wise.

Moreover, the paths of $\Sigma$ are c\`adl\`ag and of finite variation
on compacts.
\end{longlist}
\end{Theorem}

Formula \eqref{eq:intSigma} is of particular interest in connection
with stochastic volatility modeling, as in this case the integrated
volatility $\Sigma_t^+$ is a quantity of fundamental importance (see
Section \ref{sec:finecon}).

Finally, we consider the existence of moments and the second-order
structure which follow immediately from Theorems \ref{th:exmom} and
\ref{th:secord}.

\begin{Proposition}\label{th:exmompsd}
Let $\Sigma$ be a stationary $\mathbb{S}_d^+$-valued supOU process driven
by a L\'evy basis $\Lambda$ satisfying the conditions of Theorem \ref
{th:exsupoupsd}.\vspace*{-6pt}
\begin{longlist}[(iii)]
\item[(i)] If
%
\begin{equation}\label{eq:excondmoms2psd}
\int_{\|x\|>1}\|x\|^r\nu(dx)<\infty
\end{equation}
for $r\in(0,1]$, then $\Sigma$ has a finite $r$th moment, that is,
$E(\|\Sigma_t\|^r)<\infty$.

\item[(ii)] If $r\in(1,\infty)$ and
%
\begin{equation}\label{eq:excondmompsd}
\int_{\|x\|>1}\|x\|^r\nu(dx)<\infty,\qquad\int_{M_d^-}\frac{\kappa
(A)^{2r}}{\rho(A)}\pi(dA)<\infty,
\end{equation}
then $\Sigma$ has a finite $r$th moment, that is, $E(\|\Sigma_t\|
^r)<\infty$.

\item[(iii)] If the conditions given in \textup{(ii)} are satisfied for $r=2$, then the
second-order structure of $\Sigma$ is given by
\begin{eqnarray}
E(\Sigma_0)&=&-\int_{M_d^-}\mathbf{A}(A)^{-1}\biggl(\gamma_0+\int
_{\mathbb{S}
_d}x\nu(dx)\biggr)\pi(dA),\nonumber
\\
\operatorname{var}(\operatorname{vec}(\Sigma_0))&=&-\int
_{M_d^-}(\mathscr{A}(A))^{-1}
\biggl(\int_{\mathbb{S}_d}\operatorname{vec}(x)\operatorname{vec}(x)^*\nu
(dx)\biggr)\pi(dA),\nonumber%
\\
\operatorname{cov}(\operatorname{vec}(\Sigma_h),\operatorname
{vec}(\Sigma_0))&=&-\int
_{M_d^-}e^{(A\otimes I_d+I_d\otimes A)h}(\mathscr{A}(A))^{-1}\nonumber
\\
&&\hspace*{28pt}{}\times \biggl(\int_{\mathbb{S}_d}\operatorname{vec}(x)\operatorname{vec}(x)^*\nu
(dx)\biggr)\pi(dA)\nonumber
\\
\eqntext{\forall
h\in\mathbb{R}^+,}
\end{eqnarray}
with $\mathbf{A}(A)\dvtx M_d(\mathbb{R})\to M_d(\mathbb{R}), X\mapsto
AX+XA^*$ and
$\mathscr{A}(A)\dvtx M_{d^2}(\mathbb{R})\to M_{d^2}(\mathbb{R}),\break X\mapsto
(A\otimes
I_d+I_d\otimes A)X+X(A^*\otimes I_d+I_d\otimes A^*)$.
\end{longlist}
\end{Proposition}

Examples \ref{ex31}--\ref{ex5} can all be immediately adapted to
the positive semi-definite set-up. More examples in connection with
stochastic volatility modeling can be found in \cite
{BarndorffetStelzer2009sv}.

\section{Areas of applications}\label{sec:appl}

In this section we discuss possible applications for our model and the
relevance of our results for them. Some of these applications are
already developed further in other work.

\subsection{Time series modeling}
In many areas of applications (e.g., telecommunication, hydrology,
economics, finance) one is confronted with time series exhibiting a
long memory behavior (see \cite{Doukhanetal2003}, for instance) or at
least a decay of the autocovariance appearing to be a polynomial decay
rather than the exponential as typically encountered in models of
Markovian nature. Moreover, often the relevant data series are
multidimensional, and multivariate models are needed in order to
understand and adequately model the dependence effects of the observed data.
For irregularly-spaced or high-frequency data as well as when intending
to look at the data at more than one frequency, it is often advisable
not to use discrete-time models (like AR(FI)MA, see, e.g., \cite
{Brockwelletal1991}), but continuous-time models. Such models can
exhibit both continuous and discontinuous sample paths. Among the ones
with continuous sample paths are continuous-time counterparts of
ARFIMA, like FICARMA (see \cite
{Brockwelletal2005,Marquardt2007,TsaiChan2005}), which also exhibit
long memory. Often it is, however, appropriate to use models with
discontinuous sample paths. In such situations it appears adequate to
use multivariate supOU processes. It should be noted that the
individual autocovariances of multivariate supOU processes do not
necessarily have to decay monotonically like in the univariate case,
but may exhibit damped sinusoidal-like and comparable behavior due to
the involved matrix exponentials. Hence, one can reproduce second-order
moment behavior which in the univariate case calls for the use of
$\operatorname{CARMA}(p,q)$-based models.

As for all such models it is a delicate issue to estimate multivariate
supOU processes from observed data, and this poses many challenging
questions. Our calculation of the second-order moment structure clearly
opens the door for (general) method-of-moments-based techniques. Let us
illustrate this in a simple example.

Consider the set-up of Example \ref{ex31} with $B$ restricted to have
only real eigenvalues. Note that in this set-up the parameters cannot
be identified, because replacing $\beta$ by $\beta c$ and $B$ by $cB$
leaves the law of the L\'evy basis invariant for any $c>0$. So we
assume $\beta=1$ without loss of generality. Recall that
\[
\operatorname{acov}(h):=\operatorname{cov}(X_h,X_0)=-\frac{
(I_d -Bh)^{1-\alpha
}}{\alpha-1}\mathscr{B}^{-1}\biggl(\Sigma+\int_{\mathbb
{R}^d}xx^*\nu
(dx)\biggr)
\]
with $\mathscr{B}\dvtx M_d(\mathbb{R})\to M_d(\mathbb{R}), X\mapsto BX+XB^*$.
Assuming that $\Sigma+\int_{\mathbb{R}^d}xx^*\nu(dx)$ is invertible we
can define
\[
\Gamma_h=\operatorname{acov}(h)\operatorname
{acov}(0)^{-1}=(I_d-Bh)^{1-\alpha}
\]
which has eigenvalues of the form $f(h)=(1-\lambda h)^{1-\alpha}$ with
$\lambda\in(-\infty,0)$. So we can obtain an estimator $\hat\alpha
$ by calculating the empirical autocovariance function from data and
fitting $f$ to the maximum (or any other) eigenvalue of $\widehat
{\operatorname{acov}}
(h)\widehat{\operatorname{acov}}(0)^{-1}$ for $h=n\Delta$ with $n\in
\mathbb{N}$ and
$\Delta>0$ being the distance between subsequent observations (assumed
to be constant) by nonlinear least squares, for instance. Thereafter,
we can get an estimator for $B$ by $\hat B=I_d-
(\widehat{\operatorname{acov}}
(h)\widehat{\operatorname{acov}}(0)^{-1})^{\hat\alpha-1}$. Now it is
straightforward to also get estimators for $\gamma+\int_{\|x\|>1}x\nu
(dx)$ and $\Sigma+\int_{\mathbb{R}^d}xx^*\nu(dx)$ from the
empirical mean
and variance. One thus obtains the first two moments of the underlying
L\'evy process. If we restrict the allowed L\'evy bases such that the
first two moments of the underlying L\'evy process identify the
parameters $\gamma,\Sigma,\nu$, we have thus obtained a procedure to
estimate all parameters of our model.

To prove consistency and asymptotic normality (or other asymptotic
distributions in the case of true long memory) of the estimators one
obviously needs to understand the (highly non-Markovian) dependence
structure of our processes and establish mixing conditions or
appropriate substitutes. We hope to address this issue in future work.

In a set-up like in Example \ref{ex34} the estimation becomes much
easier because all parameters except the ones describing the dependence
between the components can be inferred from the univariate marginal
distributions and the univariate moment structure. In particular, the
parameters $\beta_i$ and $\alpha_i$ can be estimated from the
autocovariances of the individual one-dimensional series; hence, one
does not have to compute eigenvalues or matrix powers as above.

Finally, it should be noted that in Theorem \ref{th:supouex} or formula
\eqref{eq:cfint} we have given formulae for the characteristic
function of multivariate supOU processes. Hence, one can also use
estimation techniques based on the empirical characteristic function.
This is obviously best done in cases where the driving L\'evy basis is
chosen such that the integrals in Theorem \ref{th:supouex} or formula
\eqref{eq:cfint} can be calculated analytically. Moreover, these
formulae can be used to calculate higher order cumulants for the
methods of moments based estimation.

\subsection{Finance and econometrics}\label{sec:finecon}

L\'evy-based stochastic volatility models are very successfully applied
in both financial mathematics and financial econometrics, since they
capture many of the stylized facts (nonconstant, stochastic volatility
exhibiting jumps, heavy tails, volatility clustering, leverage
effect$,\ldots;$ see, e.g., \cite{Contetal2004,Guillaumeetal1997}) of
financial returns very well. One model often employed is the
Ornstein--Uhlenbeck-type stochastic volatility model introduced in
\cite{Barndorffetal2001c}. The use of supOU processes as the
volatility process, that is, the process modelling the instantaneous
(co)variance, allows one to introduce also long memory, which is
another important stylized fact, but not covered by most models, into
the model. Let us illustrate this in a simple set-up. Let $\Sigma$ be
a positive semi-definite supOU process as introduced in Section \ref
{sec:4} and satisfying the conditions of Theorem \ref{th:pathproppd}(ii) and (iii). Then the log-returns of $d$ financial assets (stocks or
currencies, for instance) are given by
%
\begin{equation}\label{svmodel}
Y_t=Y_0+\int_0^t (\mu+\Sigma_s\beta)\,ds+\int_0^t\Sigma
_s^{1/2}\,dW_s+\rho \,dL_t
\end{equation}
with $\mu,\beta\in\mathbb{R}^d$, initial log prices $Y_0$
independent of
$\Lambda$, $\rho\dvtx \mathbb{S}_d\to\mathbb{R}^d$ a linear operator
and $L$ being
the underlying L\'evy process. It should be noted that in this model
jumps in the price and the volatility always occur together which is
reasonable for financial data (see \cite{JacodTodorov2009}).

An extension of the above model has been investigated in-depth in \cite
{BarndorffetStelzer2009sv} where it is, in particular, shown that long
memory in $\Sigma$ causes long-range dependence in $Y$ and explicit
formulae for the moment structure of the (squared) returns are
obtained. This allows an in-depth econometric analysis and estimation
of the supOU stochastic volatility model comparable to what has been done in
\cite{PigorschetStelzer2006} for the multivariate OU-type stochastic
volatility model, where it was shown that that model can be estimated
and fits well to observed data, both from the stock and foreign
exchange markets.

In financial mathematics one is often interested in calculating prices
of derivatives from a given model and in determining the parameters
from option prices observed on the markets, referred to as calibration.
In most reasonably realistic models, unlike the Black--Scholes model,
one cannot obtain closed-form formulae for the prices of derivatives.
However, calculating prices via Monte Carlo simulations is typically
too time consuming. Whenever possible, better techniques to calculate
derivative prices, that is, conditional expectations of future payoffs,
are called for. One technique which proved to be very adequate in many
situations is the calculation of the prices by inverting the Laplace
transform (see \cite
{CarrMadan1999,EberleinGlauPapapantoleon2009,NicolatoVenardos2003}).
For the multivariate OU model this technique is successfully applied
in \cite{MuhleKarbePfaffelStelzer2009}.

Also, in the above given model \eqref{svmodel} one can calculate the
conditional Fourier transform of the future prices. Since some complex
values will arise, one has to be careful with the scalar products. We
use the scalar product $\langle x_1,x_2\rangle =x_2^*x_1$ on $\mathbb{R}^d$ and
$\langle X_1,X_2\rangle =\operatorname{tr}(X_2^*X_1)$ on $\mathbb{S}_d$ denoting by
$x^*$ the Hermitian
and by $\operatorname{tr}$ the trace of a matrix. For a linear
operator $^*$ also
denotes the adjoint operator in the following. Moreover, we assume that
$E(\exp(i\operatorname{tr}(\Lambda(B)^*u))=\exp
(\varphi
_\Lambda(u)\Pi(B))$ for all $u\in\mathbb{S}_d$ with $\varphi$
being the cumulant transform of the underlying L\'evy process and $\Pi
=\pi\times\Lambda$. Actually, $E(\exp(i\operatorname
{tr}(\Lambda
(B)^*u))$ exists for all $u\in M_d(\mathbb{R})+i\mathbb
{S}_d^+$ and
$\varphi$ can be extended to this domain as well.
Then it follows by similar arguments as in \cite
{MuhleKarbePfaffelStelzer2009,PigorschetStelzer2006} that
\begin{eqnarray}\label{eq:condft}
 E(e^{i Y_t^*u}|Y_0)&=&\exp\biggl\{i(Y_0+\mu t)^*u\nonumber
\\[2pt]
&&{}\hspace*{14pt}+\int_{M_d^-}\int_{-\infty}^t\varphi_\Lambda
\biggl[e^{A^*(t-s)}\biggl(\mathbf{A}(A)^ {-*}\biggl(u\beta^*+\frac
{i}{2}uu^*\biggr)\biggr)e^{A(t-s)}\nonumber
\\[-8pt]\\[-8pt]
&&{}\hspace*{14pt}-1_{(-\infty,0]}(s)e^{-A^*s}\biggl(\mathbf{A}(A)^
{-*}\biggl(u\beta^*+\frac{i}{2}uu^*\biggr)\biggr)e^{-As}
\nonumber
\\[-2pt]
&&{}\hspace*{26pt}-1_{(0,t]}(s)\biggl(\mathbf{A}(A)^ {-*}\biggl(u\beta
^*+\frac{i}{2}uu^*\biggr)-\rho^*u\biggr)\biggr]\,ds\,\pi(dA)\biggr\}
\nonumber
\end{eqnarray}
for all $u\in\mathbb{R}^d$ and $t\in(0,\infty)$. Here $\mathbf{A}(A)^
{-*}:=(\mathbf{A}(A)^ {-1})^*$ is the linear operator on
$M_d(\mathbb{R})$ given by $X\mapsto A^*X+XA$. If one restates the above
formula by representing $\varphi_\Lambda$ in terms of the L\'
evy--Khintchine triplet, one can calculate some of the integrals with
respect to $ds$ in the drift and Brownian covariance matrix part
explicitely, since $\int_{-\infty}^te^{A^*(t-s)}(\mathbf
{A}(A)^ {-*}X)e^{A(t-s)}\,ds=X$, for example. However, since this
results in rather lengthy formulae, especially in the part coming from
the L\'evy measure, we refrain from giving further details.
Note that we only condition on $Y_0$, since $\Sigma_0$ is highly
non-Markovian and, hence, not informative regarding future values of
$\Sigma$, and that equation \eqref{eq:intSigma} is essential to
obtain \eqref{eq:condft}.

However, it seems to be an important question what to condition upon in
such a non-Markovian setting and to the best of our knowledge this
issue has not been addressed so far. In \eqref{eq:condft} we basically
assume that we only know the current price. One could also assume that
one knows all historic prices and, hence, the historic values of
$\Sigma$, since they are given by the continuous quadratic variation
of the prices. Unfortunately, it seems extremely hard to understand
what happens if one conditions on all these historic prices. Another
point of view would be to say that $W$ and $\Lambda$ resemble the
information arriving at the markets and that market participants
observe all this information precisely. In that case it is appropriate
to condition upon the $\sigma$-algebra $\mathscr{G}_0$ generated by
$\Lambda$ up to time zero and $Y_0$ (which now is only assumed to be
independent of future values of $\Lambda$ and $W$) and one
obtains
\begin{eqnarray}
\label{eq:condft2}
 &&\hspace*{-5pt}E(e^{i Y_t^*u}|\mathscr{G}_0)\nonumber
\\[2pt]
&&\hspace*{-5pt}\qquad =\exp\biggl\{i
\biggl[(Y_0+\mu t)^*u\nonumber
\\[2pt]
&&{}\hspace*{-5pt}\hspace*{29pt}\qquad\quad +\operatorname{tr}
\biggl(\int_{M_d^-}\int_{-\infty}^0\mathbf{A}(A)^ {-1}
\bigl(e^{A(t-s)}\Lambda(dA,ds)e^{A^*(t-s)}\nonumber
\\[-8pt]\\[-8pt]
&&{}\hspace*{70pt}\qquad\quad\hspace*{53pt}-e^{-As}\Lambda
(dA,ds)e^{-A^*s}\bigr)
 \biggl(u\beta^*+\frac{i}{2}uu^*\biggr) \biggr)\biggr]\nonumber
\\[2pt]
&&{}\hspace*{-5pt}\qquad\quad\hspace*{18pt}+\int_{M_d^-}\int_{0}^t\varphi_\Lambda
\biggl[e^{A^*(t-s)}\biggl(\mathbf{A}(A)^ {-*}\biggl(u\beta^*+\frac
{i}{2}uu^*\biggr)\biggr)e^{A(t-s)}\nonumber
\\[2pt]
&&{}\hspace*{-5pt}\hspace*{20pt}\qquad\quad\hspace*{84pt}-\biggl(\mathbf{A}(A)^ {-*}\biggl(u\beta
^*+\frac{i}{2}uu^*\biggr)-\rho^*u\biggr)\biggr]\,ds\,\pi(dA)\biggr\}
\nonumber
\end{eqnarray}
 for all $u\in\mathbb{R}^d$ and $t\in(0,\infty)$.

It is clear that under appropriate technical conditions the conditional
Laplace transform given $Y_0$ or $\mathscr{G}_0$, respectively, exists
in a neighborhood of zero and \eqref{eq:condft} or \eqref{eq:condft2}
can be extended to hold on this neighborhood. Then one can use Laplace
transform techniques to calculate prices of financial derivatives. Like
for the standard OU-type model in \cite{MuhleKarbePfaffelStelzer2009}
specifications are called for under which some of the integrals can be
calculated explicitly, since otherwise the numerical integration takes
too long to make pricing and especially calibration feasible in
reasonable time. Observe that in \eqref{eq:condft2} one would set
\[
Z_t:=\int_{M_d^-}\int_{-\infty}^0\mathbf{A}(A)^ {-1}
\bigl(e^{A(t-s)}\Lambda(dA,ds)e^{A^*(t-s)}-e^{-As}\Lambda
(dA,ds)e^{-A^*s}\bigr)
\]
and determine $Z_t$ also by calibration to option prices. Obviously,
one can do this only for derivatives with a fixed maturity $t\in
\mathbb{R}
^+$, unless one increases the number of parameters one calibrates.

\subsection{Multivariate $\operatorname{supCAR}(\mathit{MA})$}

Ornstein--Uhlenbeck-type processes are a special case of the so-called
(multivariate) continuous time autoregressive moving-average (CARMA)
processes (see \cite{Brock1,Brock2,Marquardtetal2005}). As the
continuous time analogue of ARMA processes, CARMA processes, are a
fundamental class of processes for time series modeling in continuous
time. A $d$-dimensional supCAR($p$) process $Y$ can be defined by
\[
Y_t=(I_d,0,\ldots,0)\int_{M_{dp}^-}\int_{-\infty
}^te^{A(t-s)}(0,\ldots,0,I_d)^\mathrm{T}\Lambda(dA,ds),
\]
where $\Lambda$ is an $\mathbb{R}^{d}$-valued L\'evy basis on
$M_{dp}^-\times\mathbb{R}$ with generating quadruple $(\gamma,\Sigma
,\nu,\pi)$ with $\pi$ concentrated on the
matrices in $M_{dp}^-$ of the form
%
\begin{equation}\label{eq:armatrix}
\pmatrix{
0 & I_d & 0 & \cdots& 0\cr
\vdots& \ddots& \ddots&\ddots&\vdots\cr
\vdots&\ddots& \ddots& \ddots&0\cr
0 &\cdots& \cdots& 0&I_d\cr
-A_p & - A_{p-1} & \cdots&-A_2&-A_1
}
\end{equation}
with appropriate $d\times d$ matrices $A_i$. Clearly, ``$\int_{-\infty
}^te^{A(t-s)}(0,\ldots,0,I_d)^\mathrm{T}\Lambda(dA,\break ds)$'' is for $A$ fixed a
$\operatorname{CAR}(p)$ process (note that this is in contrast to \cite
{MarquardtJames2007} who define supCARMA processes differently), so it
is appropriate to call this process $\operatorname{supCAR}(p)$. Obviously many
properties for $Y$ follow from our results immediately, since it is
basically given by the first $d$-coordinates of a high-dimensional
supOU process, and this definition gives a possibility to extend CAR
processes allowing for long memory and jumps.

Using our techniques one can define supCARMA($p,q$) process with $q<p$
using an $\mathbb{R}^{d}$-valued L\'evy basis on $M_{dp}^-\times
M_{dp,d}\times\mathbb{R}$ with generating\vspace*{1pt} quadruple $(\gamma,\Sigma
,\nu
,\pi)$. To obtain proper supCARMA processes one demands $\pi(\mathscr
{A}\times M_{dp,d})=1$, denoting the set of matrices of the form \eqref
{eq:armatrix} by $\mathscr{A}$, and one sets
\[
Y_t=(I_d,0,\ldots,0)\int_{M_{dp}^-}\int_{-\infty
}^te^{A(t-s)}B\Lambda(dA,dB,ds).
\]
Adapting and extending our arguments one can easily obtain results for
this class of processes. For the interpretation as CARMA processes note
that the moving average coefficients have to be calculated from the
$d\times d$ blocks of $B$ by inverting the formulae given in \cite
{Marquardtetal2005}, Theorem 3.12.

\section{Conclusion}
In this paper we introduced multivariate supOU processes and obtained
various important properties of them. Furthermore, some areas of
application have been outlined and we are currently considering their
use in stochastic volatility modelling beginning in \cite
{BarndorffetStelzer2009sv}. However, there are still many important
issues concerning the supOU processes themselves which we hope to
address in future work. Of particular interest is, for example, the
development of good estimators for supOU models and to show properties
like consistency and asymptotic normality for them. This is related to
understanding better the dependence structure of supOU processes, which
are clearly not Markovian.

Likewise, we have shown that supOU processes allow to model long memory
effects (in a specific sense). However, a detailed theory of
multivariate long-range dependence needs to be developed.

\section*{Acknowledgments}
The authors are grateful to the Editor, Ed Waymire, and two anonymous
referees for helpful comments which considerably improved the paper.
This work was initiated during a visit of the authors to the Oxford-Man
Institute at the University of Oxford in December 2007; the authors are
very grateful for the hospitality and support given.

%

\printaddresses

\end{document}